\documentclass[11pt,amsfonts]{article}
\usepackage{dsfont}
\usepackage{graphicx}
\usepackage{latexsym}
\usepackage{amssymb}
\usepackage{amsmath}
\usepackage{enumerate}
\usepackage{stmaryrd}
\usepackage{mathrsfs}
\usepackage{layout}
\usepackage{cases}
\usepackage{eufrak}
\usepackage{euscript}
\newtheorem{prop}{Proposition}[section]
\newtheorem{lemma}{Lemma}[section]
\newtheorem{definition}{Definition}[section]

\newtheorem{theorem}{Theorem}[section]
\newtheorem{remark}{Remark}[section]
\newtheorem{example}{Example}[section]

\newtheorem{assumption}{Assumption}[section]

\numberwithin{equation}{section}

\newcounter{alphasect}
\def\alphainsection{0}

\let\oldsection=\section
\def\section{%
  \ifnum\alphainsection=1%
    \addtocounter{alphasect}{1}
  \fi%
\oldsection}%

\renewcommand\thesection{%
  \ifnum\alphainsection=1%
    \Alph{alphasect}%
  \else%
    \arabic{section}%
  \fi%
}%

\newenvironment{alphasection}{%
  \ifnum\alphainsection=1%
    \errhelp={Let other blocks end at the beginning of the next block.}
    \errmessage{Nested Alpha section not allowed}
  \fi%
  \setcounter{alphasect}{0}
  \def\alphainsection{1}
}{%
  \setcounter{alphasect}{0}
  \def\alphainsection{0}
}%

\newcommand{\stably}{\overset{st}{\rightarrow}}

\newcommand\independent{\protect\mathpalette{\protect\independenT}{\perp}}
\def\independenT#1#2{\mathrel{\rlap{$#1#2$}\mkern2mu{#1#2}}}
\newcommand{\E}{\mathbb{E}}

\def\real{{\mathord{{\rm I\kern-2.8pt R}}}}        
\def\inte{{\mathord{{\rm I\kern-2.8pt N}}}}

\def\sZZ{{\rm Z\kern-2.8ptem{}Z}}

\def\z{{\mathchoice
  {\sZZ}
  {\sZZ}
  {\rm Z\kern-0.30em{}Z}
  {\rm Z\kern-0.25em{}Z} }}
\def\sQQ{{\kern 0.27em \vrule height1.45ex width0.03em depth0em
          \kern-0.30em \rm Q}}
\def\qu{{\mathchoice
    {\sQQ}
    {\sQQ}
  {\kern 0.225em \vrule height1.05ex width0.025em depth0em \kern-0.25em \rm Q}
  {\kern 0.180em \vrule height0.78ex width0.020em depth0em \kern-0.20em \rm Q}
        }}
\def\sCC{{\kern 0.27em \vrule height1.45ex width0.03em depth0em
          \kern-0.30em \rm C}}
\def\complex{{\mathchoice
    {\sCC}
    {\sCC}
  {\kern 0.225em \vrule height1.05ex width0.025em depth0em \kern-0.25em \rm C}
  {\kern 0.180em \vrule height0.78ex width0.020em depth0em \kern-0.20em \rm C}
        }}

\newcommand{\R}{\mathbb{R}}

\newcommand{\ba}{\begin{array}}
\newcommand{\ea}{\end{array}}
\newcommand{\be}{\begin{equation}}
\newcommand{\ee}{\end{equation}}
\newcommand{\bea}{\begin{eqnarray}}
\newcommand{\eea}{\end{eqnarray}}
\newcommand{\beaa}{\begin{eqnarray*}}
\newcommand{\eeaa}{\end{eqnarray*}}

\def\z{\zeta}

\def\m{\mu}

\font\tenmath=msbm10 \font\sevenmath=msbm7 \font\fivemath=msbm5
\newfam\mathfam \textfont\mathfam=\tenmath
\scriptfont\mathfam=\sevenmath \scriptscriptfont\mathfam=\fivemath

\def \={{\buildrel {\rm (law)} \over =}}

\def \R{\mathbb{R}}

\def\qed{ \hfill \vrule width.25cm height.25cm depth0cm\smallskip}

\newcommand{\basa}{\begin{assumption}}
\newcommand{\easa}{\end{assumption}}

\newcommand{\bas}{\begin{assum}}
\newcommand{\eas}{\end{assum}}

\newcommand{\ignore}[1]{}
\usepackage{tikz}

\begin{document}
\renewcommand{\thefootnote}{\fnsymbol{footnote}}

\renewcommand{\thefootnote}{\fnsymbol{footnote}}
\begin{center}
{\Large{\bf Portmanteau inequalities on the Poisson space: }} \\ 
{\Large {\bf mixed regimes and multidimensional clustering}}

\medskip

\normalsize
by Solesne Bourguin\footnote{Universit\'e du Luxembourg. Facult\'e des Sciences, de la Technologie et de la Communication: Unit\'e de Recherche en Math\'ematiques. 6, rue Richard Coudenhove-Kalergi, L-1359 Luxembourg. Email: {\tt solesne.bourguin@gmail.com}} and Giovanni Peccati\footnote{Universit\'e du Luxembourg. Facult\'e des Sciences, de la Technologie et de la Communication: Unit\'e de Recherche en Math\'ematiques. 6, rue Richard Coudenhove-Kalergi, L-1359 Luxembourg. Email: {\tt giovanni.peccati@gmail.com}} \\ {\it  Universit\'e du Luxembourg}\\~\\

\end{center}

{\small \noindent {\bf Abstract}: Using Malliavin operators together with an interpolation technique inspired by  Arratia, Goldstein and Gordon (1989), we prove a new inequality on the Poisson space, allowing one to measure the distance between the laws of a general random vector, and of a target random element composed of Gaussian and Poisson random variables. Several consequences are deduced from this result, in particular: (1) new abstract criteria for multidimensional stable convergence on the Poisson space, (2) a class of mixed limit theorems, involving both Poisson and Gaussian limits, (3) criteria for the asymptotic independence of $U$-statistics following Gaussian and Poisson asymptotic regimes. Our results generalize and unify several previous findings in the field. We provide an application to joint sub-graph counting in random geometric graphs.    \\


    
\section{Introduction and framework}

\subsection{Overview}


\medskip

\noindent The aim of this paper is to prove and apply a new probabilistic inequality, involving vectors of random variables that are functionals of a Poisson measure defined on a general abstract space. This estimate -- which is formally stated in formula (\ref{e:pmintro}) below -- is expressed in terms of Malliavin operators, and basically allows one to measure the distance between the laws of a general random element and of a random vector whose components are in part Gaussian and in part Poisson random variables. 

\medskip

\noindent As we shall abundantly illustrate in the sequel, the inequality (\ref{e:pmintro}) is a genuine `portmanteau statement' -- in the sense that it can be used to directly deduce a number of disparate results about the convergence of random variables defined on a Poisson space, as well as to recover known ones. These results span a wide spectrum of asymptotic behaviors that are dealt with for the first time in a completely unified way. Apart from Malliavin calculus (that we apply in a form analogous to the one developed by Nualart and Vives in \cite{nuaviv}), our techniques involve the use of the Chen-Stein method (see e.g. \cite{AGG}), and provide a substantial refinement of several recent contributions concerning Central Limit Theorems (CLTs) and Poisson approximation results on the Poisson space (see \cite{LacPec_a, LacPec_b,Pec12, PSTU10, PZ10, lesmathias, Schulte12}). One of our main technical tools is an interpolation technique used in \cite{AGG} for proving multidimensional Poisson results. See e.g. \cite{NP07, np-book} for a discussion of the use of Stein-Malliavin techniques on a Gaussian space.

\medskip

\noindent As the title indicates, the two new main theoretical applications developed in the sequel are the following:

\begin{itemize}

\item[--] {\it Mixed limits:} Our results allow one to prove quantitative limit theorems (that is, limit theorems with explicit information on the rate of convergence), where the target distribution is a multidimensional combination of independent Gaussian and Poisson components.  This new class of approximation results is described in Section \ref{ss:pmixed}. They will be applied both to characterize the asymptotic independence of general $U$-statistics (see Section \ref{ss:indu}), and to subgraph counting in stochastic geometry (see Section \ref{ss:introrg}). By virtue of an approximation argument borrowed from \cite{DyMa}, part of the results discussed in Section \ref{ss:indu} extends to de-poissonized $U$-statistics. We will see that these general findings are particularly appealing for applications: in particular, we will prove a simple and surprising criterion (see Proposition \ref{prophelperalpha3}), allowing one to deduce asymptotic independence in many situations in which two regular Poisson functionals separately converge to a normal and a Poisson distribution.

\item[--] {\it Multi-dimensional Poisson convergence:} A particular choice of parameters in our main estimates allows one to deduce multidimensional Poisson approximation results, having moreover a {\it stable} nature -- in the classic sense of \cite{AlEag, ren}. This generalizes the one-dimensional findings of \cite{Pec12}. See Section \ref{ss:1} and Section \ref{ss:2}, respectively, for general statements and for applications to sequences of multiple Wiener-It\^o integrals, as well as for several comparisons with the CLTs established in \cite{PSTU10, PZ10}. Characterizing the convergence in distribution of random variables having a chaotic nature (both in a classic and a free setting) has recently become a relevant direction of research (see e.g. \cite{np-book} for an overview of the many available results in a Gaussian setting, or \cite{dnn, KNPS} for several free counterparts\footnote{A complete list of the papers related to this topic can be retrieved from the constantly updated webpage \texttt{http://www.iecn.u-nancy.fr/\~{}nourdin/steinmalliavin.htm}}), and our analysis provides substantial new contributions in the case of random variables belonging to the Poisson Wiener chaos. One should also note that Poisson approximation results based on Malliavin operators have found a number of applications in stochastic geometry, see \cite{SchTh2012}.

\end{itemize}

\medskip

%

\noindent We will illustrate our findings by completely developing an application to {\it random geometric graphs}, as described in Section \ref{ss:introrg} and Section \ref{s:fullrg}. In particular, two results will be achieved:
\begin{itemize}

\item[ (i)] a new bound for the multidimensional Poisson approximation of subgraph-counting statistics;  

\item[(ii)] a proof of a new mixed limit theorem involving the joint convergence of vectors of subgraph-counting statistics exhibiting both a Poisson and a Gaussian behavior. 
 
 \end{itemize}
 Our results extend several findings in the field -- see \cite{BhGh, JJ, Pen03}.  

\begin{remark}{\rm We anticipate some remarks about our choices in the forthcoming applications.
\begin{itemize}

\item[--] It is worth pointing out that, albeit Theorem \ref{t:indu} on the asymptotic independence of $U$--statistics is only stated for a two--dimensional vector composed of a Gaussian part and a Poisson part, it could also be stated for vectors of any dimension (both for the Gaussian and the Poisson parts) but its proof would become significantly more technical. We choose not to pursue such a level of generality, in order to keep the length of the paper within reasonable bounds.

\item[--] Similarly, Theorem \ref{t:mainrg} is only stated for a vector containing a one--dimensional Gaussian component. There is no additional difficulty in considering the higher dimensional case (using for instance the results from \cite[Chapter 3]{Pen03}) other than significantly increasing the length and technicality of the proof.

\item[--] Even though it is in principle possible to deduce the convergence results in Theorem \ref{t:mainrg} (about random graphs) by directly using Theorem \ref{t:indu}, this approach would lead to worse rates of convergence: for this reason, we prefer to provide a direct proof of all results concerning random graphs.
 
 \end{itemize}
}\end{remark}


\noindent The remainder of the paper is organized as follows. The next subsection contains a formal description of our framework: it is mostly standard material, so that {\it someone already familiar with the notation of {\rm \cite{LacPec_a, LacPec_b,PSTU10, PZ10}} can skip it at first reading}. Section \ref{s:diss} contains a detailed discussion of the main theoretical results of the paper, as well as of the applications. Section \ref{s:general} is devoted to the proofs of our general theorems, whereas Section \ref{s:fullrg} contains the proofs of our results about random graphs. An Appendix contains basic notions about Malliavin operators and contractions.

\subsection{Framework}

\noindent In what follows, we shall denote by $(Z,\mathscr{Z},\mu) $ a measure space such that
$Z$ is a Borel space, $\mathscr{Z}$ is the associated Borel $\sigma$-field, and $\mu$ is a $\sigma$-finite
Borel measure with no atoms. We write $\mathscr{Z}_{\mu} = \{ B\in \mathcal{Z}: \mu(B)< \infty \}$. The notation
$\eta = \{\eta(B) : B\in \mathscr{Z}_{\mu} \} $ is used to indicate a {\it Poisson measure} on $(Z,\mathcal{Z}) $ with {\it control} (or {\it intensity}) $\mu$. This means that $\eta $ is a collection of random variables defined on some probability space $(\Omega, \mathscr{F}, \mathbb{P}) $, indexed by
the elements of $\mathscr{Z}_{\mu} $ and such that: (i) for every $B,C \in \mathcal{Z}_{\mu}$ such that $B \cap C = \varnothing$, the random variables $ \eta(B)$ and $ \eta(C)$ are independent;  (ii) for every $B \in \mathscr{Z}_{\mu} $, $\eta(B)$ has a Poisson distribution with mean $\mu(B)$. We shall also write $$\hat{\eta}(B) = \eta(B) - \mu(B), \quad B\in \mathscr{Z}_\mu,$$ and $\hat{\eta} = \{\hat{\eta}(B) : B\in \mathscr{Z}_{\mu} \}$. A random measure satisfying property (i) is usually called ``completely random'' or ``independently scattered'' (see e.g. \cite{PeTa, SchWei} for a general introduction to these concepts, and for a discussion of any unexplained definition or result).  

\begin{remark}[The probability space] \label{rmk1}
\rm{ 

\begin{itemize}
\item[\rm (i)] In view of the assumptions on the space $(Z,\mathscr{Z},\mu)$, and to simplify the discussion,
we will assume throughout the paper that $(\Omega,\mathscr{F},\mathbb{P})$ and $\eta$ are such that
\[ \Omega = \left\{ \omega = \sum_{j=1}^{n} \delta_{z_j},n\in   \mathbb{N} \cup \{\infty\},z_j\in Z    \right\}, \]
where $\delta_z$ denotes the Dirac mass at $z$, and $\eta$ is defined as the
\textit{canonical mapping}
\[(\omega,B) \mapsto \eta(B)(\omega) = \omega(B) ,\quad B\in \mathscr{Z}_{\mu},\quad \omega\in\Omega.\]  Also, the $\sigma$-field $\mathscr{F}$ will be always supposed to be the $\sigma$-field generated by $\eta$, and we will write $L^2(\mathbb{P}) = L^2(\Omega, \mathscr{F},\mathbb{P})$. Note that the fact that $\mu$ is non-atomic implies that, for every $x\in Z$, $\mathbb{P}\{\eta\{x\} =0 \text{  or  }  1\} = 1$ .

\item[\rm (ii)] As usual, by a slight abuse of notation, we shall often write $ x \in \eta$ in order to indicate that the point $x\in Z$ is charged by the random measure $\eta(\cdot)$.
\end{itemize}

}
\end{remark}
Throughout the paper, for $p\in [1,\infty)$, the symbol $L^p(\mu)$ is shorthand for $L^p(Z,\mathscr{Z},\mu)$.
  For an integer $q\geq 2$, we shall write $L^p(\mu^q) := L^p(Z^q, \mathscr{Z}^{\otimes q}, \mu^{q}) $, whereas  $L^p_s(\mu^q)$ stands for the subspace of $L^p(\mu^q)$ composed of functions that are $\mu^{q}$-almost everywhere symmetric. Also, we adopt the convention $L^p(\mu) = L_s^p(\mu) =L^p(\mu^1) =L_s^p(\mu^1) $ and use the following standard notation: for every $q\geq 1$ and every $f,g\in L^2(\mu^q)$,
  $$ \langle f,g \rangle_{L^2(\mu^q)} = \int_{Z^q} f(z_1,...,z_q)g(z_1,...,z_q)\mu^q (dz_1,...,dz_q), \quad \|f\|_{L^2(\mu^q)} = \langle f,f \rangle^{1/2}_{L^2(\mu^q)} . $$
For every $f\in L^2(\mu^q)$, we denote by $\widetilde{f}$ the canonical symmetrization of $f$, that is,
      \[\widetilde{f}(x_1,\ldots,x_{q})=\cfrac{1}{q!}
\sum_\sigma f (x_{\sigma(1)},\ldots,x_{\sigma(q)}), \]
where $\sigma $ runs over the $q! $ permutations of the
set $\{1,\ldots,q \}$. Note that $
\|\tilde{f}\|_{L^2(\mu^q)} \leq  \|f\|_{L^2(\mu^q)}$ (to see this, use for instance the triangular inequality) .


\begin{definition}{\rm
For every deterministic function $h\in L^2(\mu)$, we write
\[I_1(h)=\hat{\eta}(h) = \int_Z h(z) \hat{\eta}(dz) \] to indicate the {\it Wiener-It\^o
integral} of $h$ with respect to $\hat{\eta}$. For every $q\geq 2$ and every $f\in L_s^2(\mu^q)$, we denote by $I_q(f)$
the {\it multiple Wiener-It\^o integral}, of order $q$, of $f$ with respect to $\hat{\eta}$. We also set $I_q(f)=I_q(\tilde{f})$, for every $f\in L^2(\mu^q)$ (not necessarily symmetric), and $I_0(b)=b$ for every real constant $b$.
}
\end{definition}
The reader is referred for instance to \cite[Chapter 5]{PeTa} or \cite{privaultbook} for a complete discussion of multiple Wiener-It\^o integrals and their properties (including the forthcoming Proposition \ref{P : MWIone} and Proposition \ref{P: MWIchaos}).

\begin{prop}\label{P : MWIone}
The following equalities hold for every $q,m\geq 1$, every $f\in L_s^2(\mu^q)$  and every $g\in L_s^2(\mu^m)$:
\begin{enumerate}
  \item[\rm 1.] $\mathbb{E}[I_q(f)]=0$,
  \item[\rm 2.] $\mathbb{E}[I_q(f) I_m(g)]= q!\langle f,g  \rangle_{L^2(\mu^q)} \mathds{1}_{\left\lbrace q=m\right\rbrace } $
  { (isometric property).}
\end{enumerate}
\end{prop}
The Hilbert space composed of the random variables of the form $I_q(f)$, where $q\geq 1$ and $f\in L^2_s(\mu^q)$, is called the $q$th \emph{Wiener chaos} associated with the Poisson measure $\eta$. The following well-known {\it chaotic representation property} is an essential feature of Poisson random measures. Recall that $\mathscr{F}$ is assumed to be generated by $\eta $.

\begin{prop}
[Wiener-It\^o chaotic decomposition] \label{P: MWIchaos} Every random variable $F\in L^2(\mathbb{P})$
admits a (unique) chaotic decomposition of the type
\begin{equation} \label{e:chaos}
F= \mathbb{E}[F] + \sum_{i = 1}^{\infty} I_i(f_i),
\end{equation}
where the series converges in $L^2(\mathbb{P})$ and, for each $i\geq 1$, the kernel $f_i$ is an element
of $L^2_s(\mu^i)$.
\end{prop}

\begin{remark}[About Malliavin calculus] {\rm For the rest of the paper, we shall use definitions and results related to Malliavin-type operators defined on the space of functionals of the Poisson measure $\eta$. Our formalism is the same as in Nualart and Vives in \cite{nuaviv}. In particular, we shall denote by 
$$
D,\, \delta,\,  L \text{\ \ and \ \ } L^{-1},
$$
 respectively, the {\it Malliavin derivative}, the {\it divergence operator}, the {\it Ornstein-Uhlenbeck generator} and its {\it pseudo-inverse}. The domains of $D$, $\delta$ and $L$ are denoted by ${\rm dom} D$, ${\rm dom} \delta$ and ${\rm dom} L$. The domain of $L^{-1}$ is given by the subclass of $L^2(\mathbb{P})$ composed of centered random variables. For the convenience of the reader we have collected some crucial definitions and results in Section \ref{app:mall} of the Appendix. Here, we just recall that, since the underlying probability space $\Omega$ is assumed to be the collection of discrete measures described in Remark \ref{rmk1}, then one can meaningfully define the random variable $\omega\mapsto F_z (\omega) =F(\omega + \delta_z),\, \omega \in \Omega, $  for every given random variable $F$ and every $z\in Z$, where $\delta_z$ is the Dirac mass at $z$. One can therefore prove the following neat representation of $D$ as a {\it difference operator} is in order: for each $F\in {\rm dom} D$,
\begin{equation}\label{e:diffop}
D_z F = F_z - F ,\,\,  \text{a.e.-} \mu(dz). 
\end{equation}
Observe that the notation $F_z(\omega) = F(\omega +\delta_z)$ extends canonically to multivariate random elements. A complete proof of this point can be found in \cite{nuaviv}. }
\end{remark}
The next statement contains an important {\it product formula} for Poisson multiple integrals (see e.g. \cite{PeTa} for a proof). Note that the statement involves contraction operators of the type $\star_r^l$: the reader is referred to Appendix \ref{app:cont} for the definition of these operators, as well as for a discussion of some relevant properties. 
\begin{prop}
[Product formula] Let $f\in L^2_s(\mu^p) $ and $g\in
L^2_s(\mu^q)$, $p,q\geq 1 $, and suppose moreover that $f \star_r^l g
\in L^2(\mu^{p+q-r-l})$ for every $r=1,\ldots,p\wedge q $ and $
l=1,\dots,r$ such that $l\neq r $. Then,
\begin{equation} \label{e:product}
I_p(f)I_q(g) = \sum_{r=0}^{p\wedge q} r!
\left(
\begin{array}{c}
  p\\
  r\\
\end{array}
\right)
 \left(
\begin{array}{c}
  q\\
  r\\
\end{array}
\right)
 \sum_{l=0}^r
 \left(
\begin{array}{c}
  r\\
  l\\
\end{array}
\right)  I_{p+q-r-l} \left(\widetilde{f\star_r^l g}\right),
\end{equation}
 with the tilde $\sim$ indicating a symmetrization.

\end{prop}

\begin{assumption}[Technical assumptions on kernels]\label{a:tech}{\rm In the sequel, whenever we consider a random vector of the type
\[
(I_{q_1}(f_1),...,I_{q_d}(f_d) ), \,\, \text{where} \,\, d\geq 1,\,\,\, q_i\geq 1,\,\,\, f_i \in L^2_s(\mu^{q_i}),
\]
we will implicitly assume that the following properties (1)-(3) are satisfied.
\begin{enumerate}

\item[(1)] For every $i=1,...,d$ and every $r=1,..., q_i$, the kernel $f_i\star_{q_i}^{q_i-r} f_i $ is an element of $L^2(\mu^{r})$.

\item[(2)] For every $i$ such that $q_i\geq 2$, every contraction of the type $(z_1,...,z_{2q_i - r- l})\mapsto |f_i|\star_r^l |f_i| (z_1,...,z_{2q_i - r- l})$ is well-defined and finite for every $r=1,...,q_i$, every $l=1,...,r$ and every $(z_1,...,z_{2q_i - r- l})\in Z^{2q_i-r-l}$.

\item[(3)] For every $ i,j=1,...,d$ such that $\max(q_i,q_j) >1$, for every $k = |q_i - q_j| \vee 1,..., q_i+q_j-2$ and every $(r,l)$ satisfying $k = q_i+q_j -2-r-l$, 
\[
\int_Z \left[\sqrt{ \int_{Z^k} (f_i(z,\cdot)\star_r^l f_j(z,\cdot))^2 \,\,d\mu^k  }\,\,\,\right]\mu(dz)<\infty,
\]
where, for every fixed $z\in Z$, the symbol $f_i(z,\cdot)$ denotes the mapping $(z_1,...,z_{q-1}) \mapsto f_i(z,z_1,...,z_{q-1})$.
\end{enumerate}

}
\end{assumption}

\begin{remark}{\rm According to \cite[Lemma 2.9 and Remark 2.10]{PZ10}, Point (1) in Assumption \ref{a:tech} implies that the following properties (a)-(c) are satisfied:

\begin{enumerate}

\item[(a)] for every $1\leq i<j\leq k$, for every $r=1,...,q_i\wedge q_j$ and every $l=1,...,r$, the contraction $f_i \star_r^l f_j$ is a well-defined element of $L^2(\mu^{q_i+q_j-r-l})$;

\item[(b)] for every $1\leq i\leq j\leq k$ and every $r=1,...,q_i$, $f_i\star_r^0 f_j$ is an element of $L^2(\mu^{q_i+q_j-r})$;

\item[(c)] for every $i=1,...,k$, for every $r=1,...,q_i$, and every $l=1,...,r\wedge (q_i-1)$, the kernel $f_i\star_r^l f_i $ is a well-defined element of $L^2(\mu^{2q_i-r-l})$.
\end{enumerate}
In particular, every random vector satisfying Assumption \ref{a:tech} is such that $I_{q_i}(f_i)^2 \in L^2(P)$ for every $i=1,...,k$,. Note that Assumption \ref{a:tech} is satisfied whenever the kernels $f_i$ are bounded functions with support in a rectangle of the type $B\times \dots\times B$, $\mu(B) <\infty$. 
}
\end{remark}

\section{Discussion of the main results}\label{s:diss}

\subsection{General bounds and mixed regimes}\label{ss:pmixed}

\noindent Fix two integers $d,m$. Observe that, in the discussion to follow, one can take either $d$ or $m$ to be zero, and in this case every expression involving such an index is set equal to zero by convention. Our main results involve the following objects:

\begin{enumerate}

\item[--] A vector ${\bf \lambda}_d = (\lambda_1,...,\lambda_d)$ of strictly positive real numbers, as well as a random vector $${\bf X}_d = (X^{(1)},...,X^{(d)}) \sim {\rm Po}_d(\lambda_1,...,\lambda_d),$$ that is, the elements of ${\bf X}_d$ are independent and such that $X^{(i)}$ has a Poisson distribution with parameter $\lambda_i$, for every $i=1,...,d$.

\item[--] A $m\times m$ covariance matrix $C = \{C(i,j) : i,j=1,...,m\}$, and a vector ${\bf N}_m = (N^{(1)},...,N^{(m)}) \sim \mathscr{N}_m (0,C)$, that is, ${\bf N}_m$ is a $m$-dimensional centered Gaussian vector with covariance $C$. We will write $H$ to indicate the $(d+m) $-dimensional random element
\begin{equation}\label{e:H}
H = ({\bf X}_d , {\bf N}_m).
\end{equation}
We shall also assume that ${\bf X}_d  \independent {\bf N}_m$, where the symbol ``$\independent$'' indicates stochastic independence, and also that $H\independent \eta$, where $\eta$  is the underlying Poisson measure.

\item[--] A vector $\mathbf {F}_d = (F^{(1)}, ....,F^{(d)})$ of random variables with values in $\mathbb{Z}_+$ such that, for every $i=1,...,d$, $F^{(i)} \in{\rm dom} D$ and $\mathbb{E}(F_i) = \lambda_i$. 

\item[--] A vector $\mathbf{G}_m = (G^{(1)},...,G^{(m)})$ of centered elements of ${\rm dom}D$. We use the notation
\begin{equation}\label{e:V}
V = ({\bf F}_d , {\bf G}_m).
\end{equation}
\noindent Note that, by definition, $V$ is $\sigma(\eta)$-measurable. 
\end{enumerate}

\begin{remark}{\rm Every asymptotic result stated in the present paper continues to hold if one allows the Poisson measure $\eta$, as well as the underlying Borel measure space $(Z, \mathscr{Z}, \mu)$, to depend on the parameter $n$ diverging to infinity.
}
\end{remark}

\noindent Our principal statement consists in an inequality allowing one to measure the distance between the laws of $H$ and $V$. To do this, we shall need the following quantities, that are defined in terms of the Malliavin operators introduced above:

\begin{eqnarray} \label{e:alpha1}
\alpha_1({\bf \lambda}_d, {\bf F}_d) &:=& \sum_{i=1}^{d}\E \left\lvert \lambda_{i} - \left\langle DF^{(i)}, -DL^{-1}F^{(i)} \right\rangle _{L^{2}(\mu)} \right\rvert \\
 \label{e:alpha2}
 \alpha_2({\bf F}_d) &:=&\sum_{i=1}^{d} \mathbb{E} \int_{Z}\left\lvert D_{z}F^{(i)}\left(D_{z}F^{(i)} - 1 \right)D_{z}L^{-1}F^{(i)} \right\rvert \mu(dz)\\
 \label{e:alpha3}
\alpha_3({\bf F}_d)&:=& \sum_{1\leq i\neq j \leq d}\E \left\lvert \left\langle DF^{(i)}, -DL^{-1}F^{(j)} \right\rangle _{L^{2}(\mu)} \right\rvert \\
\notag &&+\sum_{1\leq i\neq j\leq d} \mathbb{E} \int_{Z}\left\lvert D_{z}F^{(j)}\left(D_{z}F^{(j)} - 1 \right)D_{z}L^{-1}F^{(i)} \right\rvert \mu(dz) \\
\notag &&+  \sum_{1\leq j\neq k \leq d} \sum_{i=1}^{d} \mathbb{E} \int_{Z}\left\lvert D_{z}F^{(j)} D_{z}F^{(k)} D_{z}L^{-1}F^{(i)} \right\rvert \mu(dz)
\\ 
\label{e:beta}
\beta({\bf F}_d,{\bf G}_m) &:=& \sum_{i=1}^{d}\sum_{j=1}^{m}\E   \left\langle  | DL^{-1}G^{(j)} | ,  | DF^{(i)} | \right\rangle _{L^{2}(\mu)}  \\
\label{e:gamma1}
\gamma_1(C, {\bf G}_m) &:=& \sum_{k,j=1}^{m}\E \left\lvert C(j,k) - \left\langle DG^{(j)}, -DL^{-1}G^{(k)} \right\rangle _{L^{2}(\mu)} \right\rvert \\
\label{e:gamma2}
\gamma_2({\bf G}_m) &:=& \E \int_{Z}\left(\sum_{j=1}^{m}\left|D_{z}G^{(j)}\right|\right)^{2}\left(\sum_{j=1}^{m} \left\lvert D_{z}L^{-1}G^{(j)}   \right\rvert \right) \mu(dz).
\end{eqnarray}
As we will illustrate in great detail below, the coefficients introduced in (\ref{e:alpha1})--(\ref{e:gamma2}) should be interpreted as follows: (i) the sum $\alpha_1({\bf \lambda}_d, {\bf F}_d)+\alpha_2({\bf F}_d)$ has the form $\sum_{i=1}^d a_i$, where each $a_i$ measures the distance between the laws of $F^{(i)}$ and $X^{(i)}$, (ii) $\alpha_3({\bf F}_d)$ measures the independence between the elements of ${\bf F}_d$, (iii) the sum $ \gamma_1(C, {\bf G}_m)+\gamma_2({\bf G}_m)$ measures the distance between the laws of ${\bf G}_m$ and ${\bf N}_m$, and (iv) $\beta({\bf F}_d,{\bf G}_m)$ provides an estimate of how independent ${\bf F}_d$ and ${\bf G}_m$ are. Observe that ${\bf \lambda}_d$ and $C$ appear, respectively, only in $\alpha_1$ and $\gamma_1$. Also, one should note the asymmetric roles played by ${\bf G}_m$ and ${\bf F}_d$ in the definition of $\beta({\bf F}_d,{\bf G}_m)$.

\begin{remark}\label{r:cov}{\rm A further connection between the quantity (\ref{e:beta}) and the `degree of independence' of ${\bf F}_d$ and ${\bf G}_m$ can be obtained by combining the integration by parts formula of Lemma \ref{L : IBP} with the standard relation $L = -\delta D $, yielding that, for every $j=1,...,m$ and $i=1,...,d$,
\[
\E \left[ \left\langle DG^{(j)},  - DL^{-1}F^{(i)} \right\rangle _{L^{2}(\mu)}\right] = \E \left[ \left\langle -DL^{-1}G^{(j)},   D F^{(i)} \right\rangle _{L^{2}(\mu)}\right] = {\rm Cov}(G^{(j)} , F^{(i)} ).
\] 
A similar remark applies to the terms in $\alpha_3({\bf F}_d)$. The fact that the dependence structure of the elements of the vector $V$ can be assessed by means of a small number of parameters is a remarkable consequence of the use of the Stein and Chen-Stein methods, as well as of the integration by parts formulae of Malliavin calculus. In general, characterizing independence on the Poisson space is a very delicate (and mostly open) issue -- see e.g. \cite{privault96, privault11, RosSam}.

}\end{remark}

\medskip
\noindent We are now ready to state the main result  of the paper, namely Theorem \ref{t:mainboundintro}. The remarkable fact pointed out in its statement is that the above introduced coefficients can be linearly combined in order to measure the overall proximity of the laws of $H$ and $V$. Observe that the estimate (\ref{e:pmintro}) involves an ``adequate'' distance $d_\star(H, V)$ between the laws of the $\R^{d+m}$-valued random elements $H$ and $V$. The exact definition of such a distance (which will be always a distance providing a stronger topology than the one of convergence in distribution on $\R^{d+m}$) depends on the values of the integers $d,m$, as well as on the nature of the covariance matrix $C$, and will be formally provided in Section \ref{s:general} (see, in particular, Definition \ref{d:d1} and Definition \ref{d:d3}). For the rest of the paper, we will use the symbol ``$\,\, \stackrel{law}{\rightarrow}\,\, $'' to indicate convergence in distribution.

\begin{theorem}[Portmanteau inequality and mixed limits]\label{t:mainboundintro} Let the above assumptions and notation prevail. 
\begin{enumerate}

\item For every $d,m$ there exists an adequate distance $d_\star(\cdot, \cdot)$, as well as a universal constant $K$ (solely depending on ${\bf \lambda}_d$ and $C$), such that
\begin{equation}\label{e:pmintro}
d_\star(H,V) \leq K\left\{ \alpha_1({\bf \lambda}_d, {\bf F}_d)+\alpha_2({\bf F}_d)+\alpha_3({\bf F}_d) + \beta({\bf F}_d,{\bf G}_m) + \gamma_1(C, {\bf G}_m)+\gamma_2({\bf G}_m)\right \}. 
\end{equation} 
\item Assume $H_n = ({\bf F}_{d,n} , {\bf G}_{m,n}) $, $n\geq 1$, is a sequence of $(d+m)$-dimensional random vectors such that: (a) for every $n$, ${\bf F}_{d,n} = (F_n^{(1)},..., F_n^{(d)})$ is a vector of $\mathbb{Z}_+$-valued elements of ${\rm dom}D$ satisfying $\lambda_i(n) :=\mathbb{E}[F_n^{(i)}]   \underset{ n \rightarrow \infty}{\longrightarrow} \lambda_i$, (b) for every $n$, ${\bf G}_{m,n} = (G_n^{(1)},...,G_n^{(m)})$ is a sequence of centered elements of ${\rm dom}D$ satisfying $C_n(i,j) :=\mathbb{E}[G_n^{(i)}G_n^{(j)}] \underset{n \rightarrow \infty}{\longrightarrow} C(i,j)$ for $i,j=1,...,m$, and (c) as $n\to \infty$, 
$$
 \alpha_1({\bf \lambda}_{d,n}, {\bf F}_{d,n})+\alpha_2({\bf F}_{d,n})+\alpha_3({\bf F}_{d,n}) + \beta({\bf F}_{d,n},{\bf G}_{m,n}) + \gamma_1(C_n, {\bf G}_{m,n})+\gamma_2({\bf G}_{m,n}) \to 0,
$$
where ${\bf \lambda}_{d,n} = (\lambda_1(n),...,\lambda_d(n))$, and $C_n = \{C_n(i,j) : i,j=1,...,n\}$. Then, $H_n \stackrel{law}{\rightarrow} V$, where the convergence takes place in the sense of the distance $d_{\star}(\cdot,\cdot)$.
\end{enumerate}
\end{theorem}

\noindent The proof of Theorem \ref{t:mainboundintro}, together with a detailed statement, is provided in Section \ref{s:mainstatements}: some direct applications of the mixed  limit theorem appearing in Part 2 of its statement are described in Sections \ref{ss:introrg} and \ref{s:fullrg}, providing applications to random geometric graphs. Observe that the rest of our paper consists indeed in a series of applications of the estimate \eqref{e:pmintro}, obtained by properly selecting ${\bf \lambda}_d,\,  C,\,  {\bf F}_d$ and ${\bf G}_m$: we will use this inequality to settle a number of open questions concerning probabilistic approximations on the Poisson space. The principal theoretical applications of Theorem \ref{t:mainboundintro} developed in the present work -- namely to multidimensional Poisson approximations and asymptotic independence --  are described in the next Sections \ref{ss:stablepoisson}-\ref{ss:indu}.

\begin{remark}{\rm Specializing (\ref{e:pmintro}) to the case $m=1, \, d=0$, one obtains the main estimate in \cite{PSTU10}, concerning normal approximations of Poisson functionals in dimension one. In the case $m\geq 2, \, d=0$, (\ref{e:pmintro}) coincides with the main inequality proved in \cite{PZ10}, where the authors studied multidimensional normal approximations on the Poisson space. Finally, the case $d=0, \, m=1$ corresponds to the one-dimensional Poisson approximation result proved in \cite{Pec12}.
}
\end{remark}

\begin{remark}[About constants]\label{r:constants}{\rm By inspection of the forthcoming proof of Theorem \ref{t:mainboundintro}, the constant $K$ appearing in formula (\ref{e:pmintro}) can be taken to be have the following structure:
\begin{itemize}

\item[--] If $m =  1$ and $d\geq 1$ (in this case, $C$ is a strictly positive constant), 
$$
K = 6+ \frac{1+2\sqrt{2\pi}}{C} + \max_{i=1,...,d}\left\{ \frac{1 - e^{-\lambda_i}}{\lambda_i}+ \frac{1 - e^{-{\lambda_i}}}{\lambda^2_i}\right\} ,  
$$
where $\max_\emptyset = 0$ by convention.
\item[--] If $d \geq 1$, $m\geq 2$, then
$$
K =   11 +  \max_{i=1,...,d}\left\{ \frac{1 - e^{-\lambda_i}}{\lambda_i}+ \frac{1 - e^{-{\lambda_i}}}{\lambda^2_i}\right\} .
$$

\item[--] If $d\geq 1$ and $ m = 0$, then
$$
K =  6 \times {\bf 1}_{d>1}+  \max_{i=1,...,d}\left\{ \frac{1 - e^{-\lambda_i}}{\lambda_i}+ \frac{1 - e^{-{\lambda_i}}}{\lambda^2_i}\right\} 
$$
(the case $d=1$ follows from \cite{Pec12}).
\end{itemize} 
The values of the constants in the remaining cases (that is, when $d$ is equal to zero) can be deduced form the main results of \cite{PSTU10, PZ10}.
}
\end{remark}

\noindent We conclude this subsection with a refinement of Theorem \ref{t:mainboundintro}-2, providing useful sufficient conditions in order to have that the mixed term $\beta({\bf F}_{d,n},{\bf G}_{m,n})$ converges to zero. 

\begin{prop}
\label{prophelperalpha3}
Assume $H_n = ({\bf F}_{d,n} , {\bf G}_{m,n}) $, $n\geq 1$, is a sequence of $(d+m)$-dimensional random vectors such that ${\bf F}_{d,n} = (F_n^{(1)},..., F_n^{(d)})$ is a vector of $\mathbb{Z}_+$-valued elements of ${\rm dom}D$ and ${\bf G}_{m,n} = (G_n^{(1)},...,G_n^{(m)})$ is a sequence of centered elements of ${\rm dom}D$. Then, the following two conditions are sufficient in order to have that $\underset{n \rightarrow \infty}{\rm lim} \beta({\bf F}_{d,n},{\bf G}_{m,n}) = 0$:
\begin{itemize}

\item[--] For every $i=1,...,d$, the sequence $n\mapsto \E\left[ \int_Z \left(D_z F_n^{(i)} \right)^2 \mu(dz)\right] $ is bounded; 

\item[--] There exists $\epsilon >1$ such that, for every $j=1,...,m$, $\underset{n \rightarrow \infty}{\rm lim} \E \left[\int_Z \left |D_z L^{-1}G_n^{(j)} \right| ^{1+\epsilon} \mu(dz) \right]= 0$;


\end{itemize}
\end{prop}
{\bf Proof.} For every $i,j$, one can apply the H\"older inequality to deduce that
\begin{equation}\label{e:uff}
\E \left\langle  | DL^{-1}G^{(j)} |,| DF^{(i)}| \right\rangle _{L^{2}(\mu)} \!\! \leq \! \E\left[ \int_Z \left|D_z F_n^{(i)} \right|^{\frac{1+ \epsilon}{\epsilon} } \!\! \mu(dz)\right]^{\frac{\epsilon}{1+ \epsilon}} \!\!\times \E\left[ \int_Z \left|D_z L^{-1}G_n^{(j)} \right|^{1+ \epsilon} \!\! \mu(dz)\right]^{\frac{1}{1+ \epsilon}},
\end{equation}
and use the fact that, since $D_z F^{(i)}$ takes values in $\mathbb{Z}$, then $\left|D_z F_n^{(i)} \right|^{\frac{1+ \epsilon}{\epsilon} } \leq \left|D_z F_n^{(i)} \right|^{2}$ for every $\epsilon >1$. 

\qed

\subsection{Stable multidimensional Poisson approximations}\label{ss:stablepoisson}

We will now discuss a class of multidimensional Poisson approximation results that are a direct consequence of Theorem \ref{t:mainboundintro}. Section \ref{ss:1} contains a general statement, whereas Section \ref{ss:2} will focus on sequences of vectors of perturbed multiple integrals. We will also establish several explicit connections with the multidimensional CLTs proved in \cite{PZ10}.

\subsubsection{General statements}\label{ss:1}

As indicated in the section title, with an additional small effort we will be able to establish limit theorems in the more general framework of {\it stable convergence}. The (classic) definition of  stable convergence, in a form equivalent to the one originally given by Renyi in \cite{ren} (see also \cite{AlEag}), is provided below.
\begin{definition}[Stable convergence]\label{d:stable}{\rm 
Fix $k\geq 1$. Let $\{X_{n}\}$ be a sequence of random variables with values in $\R^{k}$, all defined on the probability space $(\Omega, \mathscr{F}, \mathbb{P})$ specified in Remark \ref{rmk1}. Let $X$ be a $\R^k$-valued random variable defined on some extended probability space $(\Omega', \mathscr{F}', \mathbb{P}')$.  We say that $X_{n}$ {\it converges  stably} to $X$,  written $X_{n} \stably X$, if 
\begin{equation*}
\underset{n \rightarrow \infty}{\rm lim}\E\left[Ze^{i\left\langle \gamma, X_{n} \right\rangle_{\R^{k}} } \right] = \E'\left[Ze^{i\left\langle \gamma, X \right\rangle_{\R^{k}} } \right] \tag{S}
\end{equation*}
for every $\gamma \in \R^{k}$ and every bounded $\mathscr{F}$--measurable random variable $Z$. }
\end{definition}

\begin{remark}\label{r:stind}{\rm In this paper, we will be exclusively interested in stable convergence results where the limiting random variable $X$ is independent of the $\sigma$-field $\mathscr{F}$. This situation corresponds to the case where $Z$ is defined on some auxiliary probability space $(A,\mathscr{A}, \mathbb{Q})$, and $(\Omega', \mathscr{F}', \mathbb{P}') = (\Omega\times A, \mathscr{F}\otimes \mathscr{A}, \mathbb{P}\otimes \mathbb{Q})$. 
}
\end{remark}
Choosing $Z=1$ in (S), one sees immediately that stable convergence implies convergence in distribution. For future reference, we now present a statement gathering together some useful results: in particular, it shows that stable convergence is an intermediate concept bridging convergence in distribution and convergence in probability. The reader is referred to \cite[Chapter 4]{JacSh} for proofs and for an exhaustive theoretical characterization of stable convergence. From now on, we will use the symbol $\stackrel{\mathbb{P}}{\rightarrow}$ to indicate convergence in probability with respect to $\mathbb{P}$.

\begin{lemma}\label{l:stable} Let $\{X_n\}$ be a sequence of random variables with values in $\R^k$.

\begin{enumerate}

\item $X_{n} \stably X$ if and only if $(X_n, Z)\stackrel{ law}{\rightarrow}(X,Z)$, for every $\mathscr{F}$-measurable random variable $Z$.

\item If $X_{n} \stably X$ and $X$ is $\mathscr{F}$-measurable, then necessarily $X_n\stackrel{\mathbb{P}}{\rightarrow}X$. 

\item If $X_{n} \stably X$ and $\{Y_n\}$ is another sequence of random elements, defined on $(\Omega, \mathscr{F}, \mathbb{P})$ and such that $Y_n \stackrel{\mathbb{P}}{\rightarrow} Y$, then $(X_n, Y_n) \stably (X,Y)$.

\item  $X_{n} \stably X$ if and only if {\rm (S)} takes place for every $Z$ belonging to a linear space $\mathcal{H}$ of bounded random variables such that $\overline{\mathcal{H}}^{L^2(\Omega, \mathscr{F}, \mathbb{P})} = L^{2}(\Omega, \mathscr{F}, \mathbb{P})$.

\end{enumerate}
\end{lemma}

\begin{remark}{\rm Properties such as Point 3 of Lemma \ref{l:stable} allow one to combine stably converging sequences with sequences converging in probability, and are one of the key tools in order to deduce limit theorems towards {\it mixtures} of probability distributions -- e.g. mixtures of Gaussian random vectors. This last feature makes indeed stable convergence extremely useful for applications, for instance within the framework of limit theorems for non-linear functionals of semimartingales, such as power variations, empirical covariances and other objects of statistical relevance. See the classic references \cite{Feigin} and \cite[Chapter 4]{JacSh}, as well as the recent survey \cite{podvetter}. Outside a semimartingale framework, stable convergence on the Wiener space has been recently studied (among others) by Peccati and Tudor in \cite{PecTud}, Peccati and Taqqu \cite{petamwii}, Nourdin and Nualart \cite{nounua} and Harnett and Nualart \cite{HarNua}. Some earlier general results about the stable convergence of non-linear functionals of random measures were obtained in \cite{PSTU08a, PeTa08b, PecTaq08c}, by using a decoupling technique known as the `principle of conditioning' -- see \cite{Jakubowki, Xue}.}
\end{remark}

\medskip

\noindent The next statement is a general stable multidimensional Poisson approximation result based on Theorem \ref{t:mainboundintro}. Recall that the {\it total variation distance} between the laws of two $\mathbb{Z}_+^d$-valued random elements ${\bf A},{\bf B}$ is given by
\begin{equation}\label{e:tv}
d_{TV} ({\bf A}, {\bf B}) = \sup_{E\subseteq \mathbb{Z}_+^d} |\mathbb{P}({\bf A}\in E) - \mathbb{P}({\bf B}\in E)|. 
\end{equation}
A proof of Theorem \ref{stablePoisson} is detailed in Section \ref{ss:poissonproof}. 
\begin{theorem}[Multidimensional stable Poisson approximations]
\label{stablePoisson}
Fix $d\geq 1$, let $(\lambda_1,...,\lambda_d)\in \R_+^d$, and let ${\bf X}_d = (X_{(1)},...,X_{(d)})\sim {\rm Po}(\lambda_1,...,\lambda_d)$ be independent of $\eta$. Let ${\bf F}_{d,n} = \left( F_{n}^{(1)},\ldots, F_{n}^{(d)} \right)$, ${n\geq 1} $, be a sequence of $\mathbb{Z}_{+}$--valued elements of ${\rm dom}\,D$ such that $\E\left[F_{n}^{(i)} \right] = \lambda_{i}(n) \underset{n \rightarrow \infty}{\rightarrow} \lambda_i$. Write ${\bf \lambda}_{d,n} = (\lambda_1(n),...,\lambda_d(n))$, $n\geq 1$, and assume moreover that: 
\begin{equation}\label{e:aaa}
\alpha_1({\bf \lambda}_{d,n}, {\bf F}_{d,n})+\alpha_2({\bf F}_{d,n})+\alpha_3({\bf F}_{d,n})  \underset{n \rightarrow \infty}{\rightarrow} 0.
\end{equation}
Then, as $n\to\infty$, the law of ${\bf F}_{d,n}$ converges to the law of ${\bf X}_{d}$ in the sense of the total variation distance, and relation (\ref{e:pmintro}) in the case $m=0$ provides an explicit estimate of the speed of convergence. If moreover,
\begin{eqnarray}\notag
\forall i = 1,\ldots, d, \, \, \forall A \in \mathscr{Z}_\mu, \, \, &&\underset{n \rightarrow \infty}{\rm lim} \E \left\lvert \int_{A} D_{z}F_{n}^{(i)} \mu(dz) \right\rvert   = \underset{n \rightarrow \infty}{\rm lim} \E \int_{A}\left\lvert D_{z}F_{n}^{(i)}(D_{z}F_{n}^{(i)}-1)  \right\rvert \mu(dz)  \\ && =  \underset{n \rightarrow \infty}{\rm lim} \E \left\lvert \int_{A} D_{z}L^{-1} F_{n}^{(i)} \mu(dz) \right\lvert =  0,\label{e:cool}
\end{eqnarray} 
and
\begin{eqnarray}\label{e:cat}
\forall 1\leq i\neq j\leq d, \, \, \forall A \in \mathscr{Z}_\mu, \, \,  \underset{n \rightarrow \infty}{\rm lim} \E \int_{A} \Big| &&\!\!\!\!\!\! \!\!\!\!\!\!D_{z}F_{n}^{(i)} D_{z}F_{n}^{(j)} \Big| \mu(dz)\\ && =  \underset{n \rightarrow \infty}{\rm lim} \E \int_{A}\left\lvert D_{z}F_{n}^{(i)}D_{z}L^{-1} F_{n}^{(i)}  \right\rvert \mu(dz)  =  0, \notag
\end{eqnarray} 
then, ${\bf F}_{d,n} \stably {\bf X}_d$.
\end{theorem}

\begin{remark}{\rm 
\begin{enumerate}
\item Theorem \ref{stablePoisson} is the first multidimensional Poisson approximation result proved by means of Malliavin operators. In the case $d=1$ (note that this implies $\alpha_3 =  0$), the fact that condition (\ref{e:aaa}) implies that $d_{TV}({\bf F}_{1,n}$, ${\bf X}_{1})\to 0$ is a consequence of the main inequality proved in \cite{Pec12}. Applications of this one-dimensional result in random geometry appear in \cite{Pec12, SchTh2012}. A new multidimensional Poisson approximation result in the context of random geometric graphs, based on the techniques developed in the present paper, appears in Theorem \ref{t:mainrg}-(c).
\item A sufficient condition (that we will be satisfied systematically in applications) in order to have that $\alpha_2({\bf F}_{d,n})+\alpha_3({\bf F}_{d,n})  \underset{n \rightarrow \infty}{\rightarrow} 0$, is that the sequences $$n\mapsto \E\left[ \int_Z \left(D_z F_n^{(i)} \right)^2 \mu(dz)\right],\quad n\mapsto \E\left[ \int_Z \left(D_z L^{-1}F_n^{(i)} \right)^2 \mu(dz)\right],  $$ are bounded for every $i$ and that, for every $i\neq j$,
$$
\underset{n \rightarrow \infty}{\rm lim} \E \int_{Z}\left\lvert D_{z}F_{n}^{(i)}(D_{z}F_{n}^{(i)}-1)  \right\rvert ^2 \mu(dz) = \underset{n \rightarrow \infty}{\rm lim} \E \int_{Z}\left\lvert D_{z}F_{n}^{(i)}D_{z}F_{n}^{(j)} \right\rvert^2 \mu(dz) =0.
$$
These conditions also imply that the middle term in (\ref{e:cool}) and the first term in (\ref{e:cat}) are equal to zero.
\item By a direct use of Point 4 of Lemma \ref{l:stable} (together with some adequate approximation argument), one can prove that another set of sufficient conditions in order to have stable convergence is that, for every $A\in \mathscr{Z}_\mu$, every $p\geq 0$ and every $i=1,...,d$,
$$
E\left [\Big| I_p(\mathds{1}_A^{\otimes p})\Big |\times  \int_{A}\left\lvert D_{z}F_{n}^{(i)}\right\rvert \mu(dz)\right]\rightarrow 0,
$$
where $\mathds{1}_A^{\otimes p} (x_1,...,x_p) = \mathds{1}_A (x_1)\cdots \mathds{1}_A (x_p)$, and $I_0 = 1$ by convention. Albeit more easily stated than (\ref{e:cool})--(\ref{e:cat}), these conditions are not simpler to verify in the applications developed in the present paper.
\end{enumerate}
}
\end{remark}

\subsubsection{The case of multiple integrals}\label{ss:2}

\noindent Now fix $d\geq 1$. Our aim is  to apply Theorem \ref{stablePoisson} in order to deduce a multidimensional Poisson approximation result for sequences of {\it perturbed multiple integrals} of the type:
\begin{equation}
\label{redriding}
{\bf F}_{d,n} =(F_n^{(1)}, ...,F_n^{(d)}) =  \left(x_{n}^{(1)} + B_{n}^{(1)} + I_{q_{1}}(f_{n}^{(1)}),\ldots, x_{n}^{(d)} + B_{n}^{(d)} + I_{q_{d}}(f_{n}^{(d)})  \right) , \quad n \geq 1,
\end{equation}
where : (i) each $F_n^{(i)}$ is a random variable with values in $\mathbb{Z}_+$, (ii) $\{x_n: n\geq 1\}$ is a sequence of positive real numbers, (iii) $q_{1},\ldots,q_{d}\geq 2$ are integers independent of $n$, (iv) $I_{q_{1}},\ldots,I_{q_{d}}$ indicate multiple Wiener-It\^o integrals of respective orders $q_{1},\ldots,q_{d}$, with respect to the compensated measure $\hat{\eta}$, (v) for each $1\leq k\leq d$, $f_{n}^{(k)} \in L_s^2(\mu^{q_{k}})$, and (vi) for each $1\leq k\leq d$, $\{B_{n}^{(k)} :n\geq 1\}$ is a {\it smooth vanishing perturbation}, in the sense of the following definition.
\begin{definition}[Smooth vanishing perturbations]\label{d:svs}{\rm A sequence $\{B_n : n\geq 1\} \subset L^2(\mathbb{P}) $ is called a {\it smooth vanishing perturbation} if $B_n, L^{-1}B_n \in {\rm dom}\, D$ for every $n\geq 1$, and the following properties hold:
\begin{eqnarray}
&& \underset{n \rightarrow \infty}{\rm lim}\E[B_n ^2] = 0\label{e:svs0}\\
\label{e:svs1}
&&\underset{n \rightarrow \infty}{\rm lim}\E\left[\|DB_n\|^2_{L^2(\mu)}\right] = \underset{n \rightarrow \infty}{\rm lim}\E\left[\|DL^{-1}B_n\|^2_{L^2(\mu)}\right] =0,\\  
&&\underset{n \rightarrow \infty}{\rm lim} \E\left[\|DB_n\|^4_{L^4(\mu)}\right] =\underset{n \rightarrow \infty}{\rm lim} \E\left[\|DL^{-1}B_n\|^4_{L^4(\mu)}\right] =  0.\label{e:svs2}
\end{eqnarray}
Note that, if (\ref{e:svs1})--(\ref{e:svs2}) are satisfied, an application of the Cauchy--Schwarz inequality yields that
\[
\underset{n \rightarrow \infty}{\rm lim}\E\left[\|DB_n\|^3_{L^3(\mu)}\right] = \underset{n \rightarrow \infty}{\rm lim}\E\left[\|DL^{-1}B_n\|^3_{L^3(\mu)}\right] =0
\]
}
\end{definition}

\begin{remark}{\rm Representing the Ornstein-Uhlenbeck semigroup as in \cite[Lemma 6.8.1]{privaultbook}, one infers that the following inequalities are always satisfied:
\[
  \E\left[\|DB_n\|^2_{L^2(\mu)}\right] \geq \E\left[\|DL^{-1}B_n\|^2_{L^2(\mu)}\right], \quad   \E\left[\|DB_n\|^4_{L^4(\mu)}\right] \geq \E\left[\|DL^{-1}B_n\|^4_{L^4(\mu)}\right]. 
\]
}
\end{remark}

\noindent The following result is the announced multidimensional Poisson approximation result for perturbed multiple integrals.

\begin{theorem}[Poisson limit theorems on perturbed chaoses] \label{perturbedStablePoisson}
Fix $d\geq 1$, $\lambda_{1},\ldots,\lambda_{d} > 0$ and let ${\bf X}_d \sim \rm{Po}_{d}(\lambda_{1},\ldots, \lambda_{d})$ be stochastically independent of $\eta$. Define the sequence ${\bf F}_{d,n}$, $n\geq 1$, according to (\ref{redriding}), and assume that for each $1 \leq i\leq d$, $x_{n}^{(i)} \underset{ n \rightarrow \infty}{\longrightarrow} \lambda_{i}$ and $\E\left[I_{q_{i}}(f_{n}^{(i)})^{2} \right] \underset{ n \rightarrow \infty}{\longrightarrow} \lambda_{i}$. Suppose also that:
\begin{equation} \label{covariance2}
 \lim_{n\rightarrow \infty} \E[F_{n}^{(i)} F_{n}^{(j)} ]=\lim_{n\rightarrow \infty}\langle f_{n}^{(i)},f_{n}^{(j)}\rangle_{L^2(\mu^{q_i})} = 0,\qquad 1\leq i\neq j \leq d.
\end{equation}
Assume moreover that the following Conditions 1--
3 hold:
\begin{enumerate}
  \item For every $k=1,...,d$, every $r=1,\ldots,q_k$ and every $l=1,\ldots, r\wedge (q_k-1)$, one has
  that $$\|f_{n}^{(k)} \star_r^l f_{n}^{(k)} \|_{L^2(\mu^{2 q_k-r-l})} \underset{ n \rightarrow \infty}{\longrightarrow} 0;$$
  \item For every $k=1,...,d$, the sequence $n \mapsto \| f_n^{(k)}\|_{L^4(\mu^{q_k})}$ is bounded and, as $n\rightarrow \infty$, $$\int_{Z^{q_{k}}} \left(  \left( f_{n}^{(k)}\right) ^2 + q_{k}!^2 \left( f_{n}^{(k)}\right) ^4 -2q_{k}! \left( f_{n}^{(k)}\right) ^3 \right) d\mu ^{q_{k}} \underset{ n \rightarrow \infty}{\longrightarrow} 0.$$
  \item For every $i\neq j$ such that $q_i = q_j$, 
  $$
  \lim_{n\rightarrow \infty}\Big \langle ( f_{n}^{(i)} )^2,(f_{n}^{(j)})^2\Big \rangle_{L^2(\mu^{q_i})} = 0.
  $$
  
\end{enumerate}
Then, ${\bf F}_{d,n} \stably {\bf X}_d$, and the convergence of ${\bf F}_{d,n}$ to ${\bf X}_d$ takes place in the sense of the total variation distance. 
\end{theorem}

\noindent A proof of Theorem \ref{perturbedStablePoisson} is provided in Section \ref{ss:poissonproof}. The following features of such a statement are noteworthy:
\begin{itemize}

\item[--] When specialized to the case $d=1$, the assumptions of Theorem \ref{perturbedStablePoisson} coincide with those in \cite[Theorem 4.1]{Pec12}. 

\item[--] In the case when $q_i \neq q_j$ for every $i\neq j$, and apart from assumption (\ref{covariance2}), the statement of Theorem \ref{perturbedStablePoisson} does not involve any requirement on the joint distribution of the elements of the vectors ${\bf F}_{d,n}$. This phenomenon mirrors some analogous findings concerning the normal approximation of vectors of multiple Wiener-It\^o integrals on the Poisson space, as first proved in \cite{PZ10}. 

\item[--] In the case where $q_i = q_j$ for $i\neq j$, Condition 3 in the statement follows automatically from (\ref{covariance2}), whenever $f_n^{(i)} $ and  $f_n^{(j)} $ have the form of a multiple of an indicator function.
\end{itemize}

\noindent For the sake of completeness, in the next statement we present a slight refinement of the chaotic CLTs proved in \cite{PZ10} (the refinement resides in the stable convergence claim). Recall that the {\it Wasserstein distance} between the laws of two $\R^m$-valued random variables $X,Y$ is given by
\begin{equation}\label{e:w1}
d_W(X,Y) = \sup_{g\in {\rm Lip}(1) }\left| \E[g(X)] - \E[g(Y)] \right|, 
\end{equation}
where ${\rm Lip}(1)$ is the class of Lipschitz functions on $\R^m$ with Lipschitz constant $\leq 1$.

\begin{theorem}[Stable CLTs for multiple integrals] \label{perturbedStableGaussian}
Fix $m\geq 1$, let ${\bf N}_m = \left( N^{(1)},\ldots, N^{(m)}\right) \sim \mathscr{N}(0,C)$, with $$C = \left\lbrace C(i,j) : i,j = 1,\ldots,m \right\rbrace $$ a $m \times m$ nonnegative definite matrix, and fix integers $q_{1},\ldots,q_{m} \geq 1$. For any $n\geq 1$ and $i = 1,\ldots,m$, let $g_{n}^{(i)} \in L_{s}^{2}(\mu^{q_{i}})$. Define the sequence ${\bf G}_{m,n} = (G_{n}^{(1)},\ldots, G_{n}^{(m)})$, $n \geq 1$, as
\[
G_n^{(i)} = I_{q_i}(g_n^{(i)}), \quad n\geq 1, \,\, i=1,...,m.
\]
Suppose that Assumption \ref{a:tech} is satisfied for every $n$, and also that
\begin{equation} \label{covariance}
 \lim_{n\rightarrow \infty} \E[G_{n}^{(i)} G_{n}^{(j)} ]=\mathds{1}_{(q_j=q_i)}\times\lim_{n\rightarrow \infty}\langle g_{n}^{(i)},g_{n}^{(j)}\rangle_{L^2(\mu^{q_i})} = C(i,j),\qquad 1\leq i,j \leq m.
\end{equation}
Assume moreover that the following Conditions 1--2 hold for every $k=1,...,m$:
\begin{enumerate}
  \item For every $r=1,\ldots,q_k$ and every $l=1,\ldots, r\wedge (q_k-1)$, one has
  that $$\|g_{n}^{(k)} \star_r^l g_{n}^{(k)} \|_{L^2(\mu^{2 q_k-r-l})} \underset{ n \rightarrow \infty}{\longrightarrow} 0;$$
  \item As $n\rightarrow \infty$, $\int_{Z^{q_{k}}}\left( g_{n}^{(k)}\right) ^4 d\mu ^q_{k} \underset{ n \rightarrow \infty}{\longrightarrow} 0$.
\end{enumerate}
Then, $\gamma_1(C_n, {\bf G}_{m,n})+\gamma_2({\bf G}_{m,n})\to 0$, ${\bf G}_{m,n} \stably {\bf N}_m$ and the convergence of ${\bf G}_{m,n}$ to ${\bf N}_m$ takes place in the sense of the Wasserstein distance.
\end{theorem}

\begin{remark}{\rm Apart from the covariance condition (\ref{covariance}), the assumptions appearing in the previous statement do not involve any requirement on the joint distribution of the components of the vector ${\bf G}_{n,m}$. Moreover, these assumptions are the same as those in \cite[Theorem 5.1]{PSTU10} (for the case $m=1$) and \cite[Theorem 5.8]{PZ10} (for the case $m\geq 2$). The somewhat remarkable (albeit easily checked) fact stated in Theorem \ref{perturbedStableGaussian} is that the same assumptions implying a CLT for multiple integrals systematically yield a stable convergence result. Note that this phenomenon represents the exact Poisson space counterpart of a finding by Peccati and Tudor \cite{PecTud}, concerning the stable convergence of vectors of multiple integrals with respect to a general Gaussian field. See \cite[Chapter 6]{np-book} for an exhaustive discussion of this phenomenon. CLTs on the Poisson space based on contraction operators have already been applied to a variety of frameworks -- such as CLTs for linear and non-linear functionals of L\'evy driven moving averages \cite{PSTU10, PecTaq08c}, characterization of hazard rates in Bayesian survival models \cite{dpp, pecpru} and limit theorems in stochastic geometry \cite{LacPec_a, LacPec_b}.  }
\end{remark}

\noindent We conclude this section by stating an application of Proposition \ref{prophelperalpha3}, implying that vectors of (perturbed) multiple integrals satisfying the assumptions of Theorem \ref{perturbedStablePoisson} and Theorem \ref{perturbedStableGaussian} are automatically independent in the limit. 

\begin{prop}\label{p:asint} Let the sequences $\{ {\bf F}_{d,n} : n\geq 1\}$ and $\{ {\bf G}_{m,n} : n\geq 1\}$, respectively, satisfy the assumptions of  Theorem \ref{perturbedStablePoisson} and Theorem \ref{perturbedStableGaussian}. Then, $\beta({\bf F}_{d,n},{\bf G}_{m,n})\to 0$, as $n\to\infty$, and the two sequences are asymptotically independent.

\end{prop}

\noindent Several connected results involving $U$-statistics are discussed in the next section.

\subsection{Asymptotic independence of $U$-statistics}\label{ss:indu}

We shall now apply the main findings of the paper in order to characterize the asymptotic independence of sequences of random variables having the form of $U$-statistics converging either to a Gaussian or a Poisson limit. Our basic message is that, under fairly general conditions, $U$-statistics satisfying a CLT are necessarily asymptotically independent of any $U$-statistic converging to Poisson. The criteria for Gaussian and Poisson convergence used below are taken from references \cite{LacPec_a, LacPec_b, lesmathias} and \cite{SchTh2012}: to our knowledge, these references contain the most general conditions in order for a sequence of $U$-statistics based on a Poisson measure to converge, respectively, to a Gaussian or a Poisson limit. 
\medskip

\noindent By virtue of a de-poissonization argument borrowed from \cite{DyMa}, we will be able to deal both with poissonized and non-poissonized $U$-statistics based on a i.i.d. sequence -- see Proposition \ref{p:zap}. The reader is referred to \cite{KoBo} for a survey of the classic theory of $U$-statistics. See \cite{BhGh,HSS,JJ,SiBr}, as well as the monograph \cite{Pen03} and the references therein, for several examples of the use of $U$-statistics in stochastic geometry. See \cite{DFR,LacPec_a,LacPec_b,Pec12,lesmathias,Schulte,Schulte12,SchTh2012} for new geometric applications based on Stein-Malliavin techniques. Albeit unified studies of Gaussian and Poisson limits for $U$-statistics are available (see e.g. \cite{JJ}), we could not find in the literature any systematic characterization of the asymptotic independence of $U$-statistics in the spirit of the present section.

\begin{remark}{\rm
Section \ref{ss:introrg} contains another characterization of asymptotic independence of $U$-statistics associated with random geometric graphs. Rather than using the general results discussed below, and due to the explicit nature of the kernels involved, we will establish such results by some direct analytical computations -- allowing to obtain better rates of convergence, as well as results in higher dimensions.
}
\end{remark}

\noindent Since it is relevant for applications, we will explicitly work with a sequence of Poisson measures $\{ \eta_n : n\geq 1\}$, each defined on the Borel space $(Z,\mathscr{Z})$ and controlled by a $\sigma$-finite measure $\mu_n$ {\it possibly depending on} $n$. Following \cite[Section 3.1]{lesmathias}, we now introduce the concept of a $U$-{\it statistic} associated with the Poisson measure $\eta_n$. 

\begin{definition} {\bf ($U$-statistics)} {\rm Fix $k\geq 1$. A random variable $F$ is called a $U$-{\it statistic of order} $k$, based on the Poisson measure $\eta_n$ with control $\mu_n$, if there exists a kernel $h\in L^1_s(\mu_n^k)$ such that
\begin{equation}\label{e:ustat}
F = \sum_{{\bf x} \in \eta^k_{n,\neq}} h({\bf x}), 
\end{equation}
where the symbol $\eta^k_{n,\neq}$ indicates the class of all $k$-dimensional vectors ${\bf x} =(x_1,\dots,x_k)$ such that $x_i\in \eta_n$ and $x_i\neq x_j$ for every $1\leq i\neq j \leq k$. As formally explained in {\rm \cite[Definition 3.1]{lesmathias}}, the possibly infinite sum appearing in (\ref{e:ustat}) must be regarded as the $L^1(\mathbb{P})$ limit of objects of the type $\sum_{{\bf x}\in \eta^k_{n,\neq}\cap A_q} f({\bf x})$, $q\geq 1$, where the sets $A_q \in Z^k$ are such that $\mu_n^k(A_q)<\infty$ and $A_q \uparrow Z^k$, as $q\to \infty$.}
\end{definition}

\begin{example}[Poissonized $U$-statistics]\label{ex:poissu}{\rm Assume $\{Y_i : i\geq 1\}$ is a sequence of i.i.d. random variables with values in $Z$ and common non-atomic distribution $p$, and consider an independent Poisson random variable $N(n)$ with parameter $n \geq 1$. Then,  $\eta_n(\cdot) = \sum_{i=1}^{N(n)}\delta_{Y_i}(\cdot)$ is a Poisson random measure with control $\mu_n = n p$. In this framework, for every $k\geq 1$ and any symmetric kernel $h \in L_s^1(\mu_n^k) = L_s^1((np)^k)$, the corresponding $U$-statistic has the form
\begin{equation}\label{e:suz}
F = \sum_{{\bf x} \in \eta^k_{n,\neq}} h({\bf x}) = \sum_{1\leq i_1,...,i_k \leq N(n); \, i_i\neq i_j } h(Y_{i_1},...,Y_{i_k}).
\end{equation}
The random variable obtained by replacing $N(n)$ with the integer $n$ in (\ref{e:suz}) is customarily called the {\it de-poissonized version} of $F$.}
\end{example}

\noindent The following crucial fact is proved by Reitzner \& Schulte in \cite[Lemma 3.5 and Theorem 3.6]{lesmathias}:

\begin{prop}
\label{p:zap}
Consider a kernel $h\in L_{s}^1(\mu_n^k)$ such that the corresponding $U$-statistic $F$ in (\ref{e:ustat}) is square-integrable. Then, $h$ is necessarily square-integrable, and $F$ admits a chaotic decomposition of the form (\ref{e:chaos}), with \begin{equation}\label{e:lastpenrose}
f_i({\bf x}_{i}) := h_i({\bf x}_{i})= \binom{k}{i}\int_{Z_n^{k-i}} h({\bf x}_{i},{\bf x}_{k-i})\, d\mu_n^{k-i},\quad {\bf x}_{i}\in Z^i,
\end{equation}
for $1 \leq i \leq k$, and $f_{i}=0$ for $i>k$. In particular, $h=f_k$ and the projection $f_{i}$ is in $ L_{s}^{1,2}(\mu_n^i)$ for each $1 \leq i \leq k$.\\
\end{prop}

\begin{remark}{\rm
In \cite{lesmathias} it is proved that the condition $h\in L^{1}(\mu_n^k)\cap L^{2}(\mu_n^k) $ does not ensure, in general, that the associated $U$-statistic $F$ in (\ref{e:ustat}) is a square-integrable random variable.} 
\end{remark}

\noindent The forthcoming Theorem \ref{t:indu} is the main result of the section. It is divided in three parts. Part 1 collects some of the main results from \cite{LacPec_a, LacPec_b} concerning the normal approximation of $U$-statistics. Part 2 contains conditions for Poisson approximations of $U$-statistics taken from \cite[Proposition 4.1]{SchTh2012}. Part 3 is new and states that, under the conditions appearing in the first two parts, any two $U$-statistics converging, respectively, to a Gaussian and a Poisson limit are necessarily asymptotically independent. 

\begin{remark}{\rm The bounds from \cite{LacPec_a, LacPec_b} stated below are easier to handle than the ones deduced in the seminal work \cite{lesmathias} -- albeit they are basically equivalent in several applications. The resulting conditions for asymptotic normality have been proved in \cite{LacPec_a} to be necessary and sufficient in many important instances. The conditions for Poisson approximations taken from \cite{SchTh2012} should be compared with the classic findings of \cite{JJ,SiBr}. 
}
\end{remark}

\noindent Our framework is the following:

\begin{itemize}
\item[--]  The sequence 
$$
G_n = \sum_{{\bf x} \in \eta^k_{n,\neq}} g_n ({\bf x}), \quad n\geq 1,
$$
is composed of square-integrable $U$-statistics of order $k\geq 2$ such that $g_n\in L^{1}(\mu_n^k)\cap L^{2}(\mu_n^k)$. We write $g_{i,n}$, $i=1,...,k$, for the $i$th kernel in the chaotic decomposition of $g_n$, as given in (\ref{e:lastpenrose}). We write $\sigma^2_n = {\rm Var} (G_n)$ and write $\tilde{ G}_n = [G_n - \E( G_n)]/\sigma_n$. 
\item[--]  For $G_n$ as above, we set
\begin{eqnarray}\label{e:B3}
&& B(G_n;\sigma_n) = \frac{1}{\sigma^2_n}\left\{ \max_{(\ast)} \ \|g_{i,n}\star_r^l g_{j,n}\|_{L^2(\mu_n^{i+j-r-l})} + \max_{i=1,...,k} \|g_{i,n}\|^2_{L^4(\mu_n^{i})}\right\},
\end{eqnarray}
where $\displaystyle{\max_{(\ast)}}$ ranges over all quadruples $(i,j,r,l)$ such that  $1 \leq l \leq r \leq i \leq j $ ($i,j \leq k$) and $l\neq j$ (in particular, quadruples such that $l =r=i=j=1$ do not appear in the argument of $\displaystyle{\max_{(\ast)}}$).

\item[--] For an integer $k'\geq 2$, $\{O_n : n\geq 1 \} $ is a sequence of symmetric elements of $\mathscr{Z}^{k'}$ such that $\mu^{k'}_n(O_n)<\infty$ for every $n$. For every $n$, we define $F_n$ to be the $U$-statistic obtained from (\ref{e:ustat}) by taking $h({\bf x}) = h_n ({\bf x})= k'!^{-1}\mathds{1}_{O_n}({\bf x})$. To simplify the discussion, we may assume that each $O_n$ is contained in a $k'$-fold Cartesian product of the type $K_n \times \cdots \times K_n$, with $\mu_n(K_n)<\infty$, thus ensuring that each $F_n$ is square-integrable. Accordingly, we denote by $h_{i,n}$, $i=1,...,k'$, the $i$th kernel in the chaotic decomposition of $F_n$, and we also write $\lambda_n = k'!^{-1}\mu_n^{k'}(O_n) = \E[F_n]$.

\item[--] Define: 
$$
\rho_n = \sup \mu_n^j\Big \{ ( y_1,...,y_j)\in Z^j : (y_1,...,y_j, o_1,...,o_{k'-j}) \in O_n\Big\}. 
$$
where the supremum runs over all $j=1,...,k'-1$ and all vectors $(o_1,...,o_{k'-j}) \in Z^{k'-j}.$.
\end{itemize}

\begin{theorem}\label{t:indu} We denote by $N$ and $X_\lambda$, respectively, a $\mathscr{N}(0,1)$ and a ${\rm Po}(\lambda)$ random variable, where $\lambda>0$. We assume that $N \independent X_\lambda$. 
\begin{itemize}

\item[1.] There exists a constant $C_k>0$, independent of $n$ such that,  
$$
d_W(\tilde{G}_n,N) \leq C_k B(G_n; \sigma_n).
$$
In particular, if $B(G_n; \sigma_n)\to 0 $, then $\tilde{G}_n$ converges in distribution to $N$, in the sense of the Wasserstein distance.
\item[2.] There exists a constant $D_{k'}>0$, independent of $n$, such that
$$
d_{TV}(F_n, X_\lambda) \leq |\lambda_n - \lambda| +D_{k'} \frac{1-e^{-\lambda_n}}{\lambda_n}\left(1+\frac{1}{\lambda_n}\right)\sqrt{(\lambda_n+\lambda_n^{2})(\rho_n +\rho_n^{4})} := A_n.
$$
In particular, if $A_n \to 0 $, then $F_n$ converges in distribution to $X_\lambda$, in the sense of the total variation distance.
\item[3.] Write $V_n = (F_n, \tilde{G}_n)$, and $H = (X_\lambda,N)$. For an adequate distance, $d_\star$, there exists a constant $M =M(k;k')$, independent of $n$ such that
$$
d_{\star} (V_n, H)\leq M\times \left\{ A_n +B(G_n; \sigma_n)+B(G_n; \sigma_n)^{1/2} \right\}
$$
In particular, if $\underset{n \rightarrow \infty}{\rm lim} A_n= \underset{n \rightarrow \infty}{\rm lim}  B(G_n; \sigma_n) = 0 $, then $V_n$ converges in distribution to $H$, and $F_n$ and $\tilde{G}_n$ are asymptotically independent.
\end{itemize}
\end{theorem}

\noindent The next statement is the announced de-Poissonization result. 

\begin{prop}[De-poissonization]\label{t:indue} Let the notation of Theorem \ref{t:indu} prevail, and assume that, for every $n$, the Poisson measure $\eta_n$ is defined as in Example \ref{ex:poissu}. Write $F^0_n$ and $\tilde{G}^0_n$ to indicate the de-poissonized versions of $F_n$ and $\tilde{G}_n$. If $\underset{n \rightarrow \infty}{\rm lim} A_n= \underset{n \rightarrow \infty}{\rm lim}  B(G_n; \sigma_n) = 0 $, then $(F^0_n, \,\tilde{G}^0_n) $ converges in distribution to $H$. 
\end{prop}

\subsection{Applications to random graphs}\label{ss:introrg}

\noindent We now demonstrate how to apply our main results to study multidimensional limit theorems for subgraph-counting statistics in the disk-graph model on $\R^m$. Our main contribution, stated in Theorem \ref{t:mainrg} below, is a new estimate providing both mixed limit theorems and multidimensional Poisson approximation results. The present section contains statements, examples and discussions; proofs are detailed in Section \ref{s:fullrg}. Our notation has been chosen in order to loosely match the one adopted in \cite[Chapter 3]{Pen03}, as well as in \cite[Section 3]{LacPec_b}.

\medskip

\noindent We fix $m\geq 1$, as well as a bounded and continuous probability density $f$ on $\R^m$. We denote by $Y = \{Y_i : i\geq 1\}$ a sequence of $\R^m$-valued i.i.d. random variables, distributed according to the density $f$. For every $n=1,2,...$, we write $N(n)$ to indicate a Poisson random variable with mean $n$, independent of $Y$. It is a standard result that the random measure $\eta_n = \sum_{i=1}^{N(n)} \delta_{Y_i}$, where $\delta_x$ indicates a Dirac mass at $x$, is a Poisson measure on $\R^m$ with control measure given by $\mu_n(dx) = n f(x)dx$ (with $dx$ indicating the Lebesgue measure on $\R^m$). We shall also write $\hat{\eta}_n = \eta_n - \mu_n$, $n\geq 1$. Given positive sequences $a_n, b_n$, we write $b_n\sim a_n$ to indicate that the ratio $a_n/b_n$ converges to 1, as $n\to\infty$.

\medskip

\noindent Let $\{t_n : n\geq 1\}$ be a sequence of strictly decreasing positive numbers such that $\underset{n \rightarrow \infty}{\rm lim} t_n = 0$. For every $n$, the symbol $G'(Y;t_n)$ indicates the undirected random {\it disk graph} obtained as follows: the vertices of $G'(Y;t_n)$ are given by the random set $V_n = \{Y_1,...,Y_{N(n)}\}$ and two vertices $Y_i,Y_j$ are connected by an edge if and only if $\| Y_i - Y_j \|_{\R^m} \in (0,t_n)$. By convention, we set $G'(Y;t_n) = \emptyset$ whenever $N(n)=0$. Now fix $k\geq 2$, and let $\Gamma $ be a connected graph of order $k$. For every $n\geq 1$, we shall denote by $G'_n(\Gamma)$ the number of {\it induced subgraphs} of $G'(Y;t_n)$ that are isomorphic to $\Gamma$, that is: $G'_n(\Gamma)$ counts the number of subsets $\{i_1,...,i_k\} \subset \{1,...,N(n)\}$ such that the restriction of $G'(Y;t_n)$ to $\{Y_{i_1},...,Y_{i_k}\}$ is isomorphic to $\Gamma$. Every graph $\Gamma$ considered in the sequel is assumed to be {\it feasible} for every $n$: this means that the probability that the restriction of $G'(Y;t_n)$ to $\{Y_1,...,Y_k\}$ is isomorphic to $\Gamma$ is strictly positive for every $n$. The study of the asymptotic behavior of the random variables $G'_n(\Gamma)$, as $n$ goes to infinity, is one of the staples of the modern theory of geometric random graphs, and many results are known. The reader is referred to Penrose \cite[Chapter 3]{Pen03} for a general discussion and for detailed proofs, and to \cite[Section 3]{LacPec_b} and \cite[Section 6]{lesmathias} for some recent refinements. 

\medskip

\noindent In what follows we will focus on the following setup: (i) $k_0,k$ are integers such that $2\leq k_0 < k$, (ii) the sequence $\{t_n \}$ introduced above is such that $t^m_n \sim  n^{-\frac{k}{k-1}}$,  (iii) $\Gamma_0$ is a feasible connected graph of order $k_0$, (iv) for some $d\geq 1$, $(\Gamma_1,...,\Gamma_d)$ is a collection of non--isomorphic feasible connected graphs with order $k$. We also write
$$
\tilde{G}'_n(\Gamma_0) = \frac{{G}'_n(\Gamma_0) - \E[{G}'_n(\Gamma_0)] }{{\rm Var}({G}'_n(\Gamma_0))^{1/2}}.
$$
The specificity of this framework is that, for such a sequence $\{t_n\}$, the random variables $\tilde{G}'_n(\Gamma_0)$ and ${G}'_n(\Gamma_j)$ ($j=1,...,d$) satisfying, respectively, a CLT and a Poisson limit theorem.  Our principal aim is to provide an exhaustive description of their joint asymptotic distribution. The following statement gathers together many results from the literature, mostly taken from \cite[Chapter 3]{Pen03} (for limit theorems, expectations and covariances) and \cite{LacPec_b} (for the estimates on the Wasserstein distance).

\begin{prop}\label{p:g} Let the above notation and assumptions prevail. 
\begin{enumerate}

\item[\rm (a) ] There exist constants $a_0,\, b_0>0$ such that, as $n\to\infty$, $$\mathbb{E}[G'_n(\Gamma_0)] \sim a_0 n^{k_0}(t_n^m)^{k_0-1} \sim a_0 n ^{(k-k_0)/(k-1)} \to \infty,$$ and ${\rm Var}({G}'_n(\Gamma_0))\sim  b_0 n ^{(k-k_0)/(k-1)} \to \infty$. Moreover, the random variable $\tilde{G}'_n(\Gamma_0)$ converges in distribution towards a $\mathscr{N}(0,1)$ random variable, with an upper bound of order $ n ^{-(k-k_0)/2(k-1)}$ on the Wasserstein distance. 

\item[\rm (b)] There exist constants $a_1,...,a_d>0$ such that, for every $j=1,...,d$ 
$$\mathbb{E}[G'_n(\Gamma_j)] \sim {\rm Var}({G}'_n(\Gamma_j) \sim a_j n^k (t_n^m)^{k-1} \to a_j.$$
Moreover, $({G}'_n(\Gamma_1),..., {G}'_n(\Gamma_d))$ converges in distribution to a $d$-dimensional vector $(X_1,...,X_d)$ composed of independent random variables such that $X_j$ has a Poisson distribution with parameter $a_j$.

\item[\rm (c) ] As $n\to \infty$, one has that, for every $i,j=1,...,d$, ${\rm Cov}(\tilde{G}'_n(\Gamma_0) , {G}'_n(\Gamma_j)) = O\left( n^{-(k-k_0)/2(k-1)}\right) $, and ${\rm Cov}({G}'_n(\Gamma_i) , {G}'_n(\Gamma_j)) = O\left( n^{-1 /(k-1)}\right) $.
\end{enumerate}
\end{prop}

\begin{remark}{\rm  We could not find a proof of the multidimensional Poisson limit theorem stated at Point (b) of the previous statement. However, such a conclusion can be easily deduced e.g. from \cite[Theorem 3.5]{Pen03}, together with a standard poissonization argument.
}
\end{remark}

\noindent Plainly, Proposition \ref{p:g} does not allow to deduce a characterization of the joint asymptotic distribution of the components of the vector $$V_n := ({G}'_n(\Gamma_1),..., {G}'_n(\Gamma_d), \tilde{G}'_n(\Gamma_0)), \quad n\geq 1.$$ In particular, albeit Part (c) of such a statement implies that the random variables $\tilde{G}'_n(\Gamma_0)$ and ${G}'_n(\Gamma_j)$ are asymptotically uncorrelated for every $j=1,...,d$, nothing can be a priori inferred about their asymptotic independence. The following statement, which provides a highly non-trivial application of Theorem \ref{t:mainboundintro}, yields an exhaustive characterization of the joint asymptotic behavior of the components of $V_n$.

\begin{theorem}[Mixed regimes in random graphs]\label{t:mainrg} For every $n$ and every $j=1,...,d$, set $\lambda_{j,n} = \E[G'_n(\Gamma_j)]$. Let $N\sim \mathscr{N}(0,1)$, ${\bf X}_{d,n} = (X_{1,n},...,X_{d,n})\sim {\rm Po}_d(\lambda_{1,n},...,\lambda_{d,n})$, assume that $N$ and ${\bf X}_{d,n}$ are stochastically independent, and write $H_n = ({\bf X}_{d,n}, N)$. 

\begin{enumerate}

\item[\rm (a)] There exist two constants $A$ and $B$, independent of $n$, such that, for some adequate distance $d_\star$,
\begin{equation}\label{e:rgbound}
d_{\star}(V_n, H_n) \leq A \sqrt{n t_n^m} + Bn^{-\frac{k-k_0}{4(k-1)}}= O\left( n^{ - \frac{1}{2(k-1)}}+n^{-\frac{k-k_0}{4(k-1)}}\right) .
\end{equation}

\item[\rm (b)] Let ${\bf X}_{d} \sim {\rm Po}_d(a_1,...,a_d)$, where the constants $a_j$ have been defined in Part (b) of Proposition \ref{p:g}, be independent of $N$, and set $H = ({\bf X}_{d}, N)$. Then, as $n\to \infty$, $V_n$ converges in distribution to $H$.

\item[\rm (c)]  Write $V'_n := ({G}'_n(\Gamma_1),..., {G}'_n(\Gamma_d)), \, n\geq 1.$ Then, there exists a constant $C$, independent of $n$, such that
\begin{equation}\label{e:rgbound2}
d_{TV}(V'_n, {\bf X}_{d,n}) \leq C \sqrt{n t_n^m} = O\left( n^{ - \frac{1}{2(k-1)}}\right) . 
\end{equation}

\end{enumerate}
\end{theorem} 

\begin{remark}{\rm 

\begin{itemize}

\item[(i)] The estimates (\ref{e:rgbound})--(\ref{e:rgbound2}) and the content of Point (b) are new. We do not know of any other available technique providing the limit theorem at Point (b). Note that such a statement yields, in particular, the asymptotic independence of $V'_n$ and $\tilde{G}'_n(\Gamma_0)$. The rate of convergence implied by formula (\ref{e:rgbound2}) is probably suboptimal (the correct rate should be of the order of $nt_n^m$ -- compare with the statement of \cite[Theorem 3.5]{Pen03} in the case of non-Poissonized graph). It is plausible that one could obtain a better rate by avoiding the use of the Cauchy-Schwarz inequality in the proof, and by estimating expectations by means of some generalized Palm-type computations (see e.g. \cite[Section 1.7]{Pen03}). This approach requires several technical computations; to keep the length of the present paper within bounds, we plan to address this issue elsewhere. Previous classic references on geometric random graphs are \cite{BhGh, JJ, SiBr}. 

\item[(ii)] A quick computation shows that if $k = k_0 +1$, the rate of convergence in (\ref{e:rgbound}) is $O\left( n^{-\frac{1}{4(k-1)}}\right) $ and if $k \geq k_0 +2$, the rate of convergence is $O\left( n^{ - \frac{1}{2(k-1)}}\right) $. 
\end{itemize}
}
\end{remark}

\begin{example}{\rm Let $k_0=2, \, k=3$, and consider the sequence of disk graphs with radius $t_n$ such that $t_n^m \sim n^{-3/2}$. Define the following graphs: (i) $\Gamma_0$ is the connected graph with two-vertices, (ii) $\Gamma_1$ is the triangle and, (iii) $\Gamma_2$ is the 3-path, that is, the connected graph with three vertices and two edges. Plainly, $G'_n(\Gamma_0)$ equals the number of edges in the disk graph, whereas $G'_n(\Gamma_1)$ and $G'_n(\Gamma_2)$ count, respectively, the number of induced triangles and of induced 3-paths. Since $\Gamma_1$ and $\Gamma_2$ are non-isomorphic, Theorem \ref{t:mainrg} can be applied, and we deduce that $\tilde{G}'_n(\Gamma_0),\, G'_n(\Gamma_1)$ and $G'_n(\Gamma_2)$ are asymptotically independent, and that they jointly converge towards a mixed Poisson/Gaussian vector, with an upper bound on the speed of convergence of the order of $n^{-1/8}$.
}
\end{example}

\noindent We conclude this section by pointing out that an application of the de-poissonization Lemma \ref{l:lalla} yields the following generalization of Theorem \ref{t:mainrg}. The details of the proof are left to the reader. For every $n$, we denote by $G(Y;t_n)$ the de-poissonized random graph obtained as follows: the vertices of $G(Y;t_n)$ are given by the random set $V_n = \{Y_1,...,Y_n\}$ and two vertices $Y_i,Y_j$ are connected by an edge if and only if $\| Y_i - Y_j \|_{\R^m} \in (0,t_n)$.

\begin{prop} The conclusion of Theorem \ref{t:mainrg}-(b) continues to hold whenever the disk graph $G'(Y;t_n)$ is replaced with the de-poissonized random graph $G(Y;t_n)$, and each counting statistic $G'_n(\Gamma_i)$ is replaced by its de-poissonized counterpart. 
\end{prop}

\section{Proofs of the main theorems}\label{s:general}

\subsection{Preliminaries}

\noindent We will now introduce several classes of functions that will be used to define particular metrics used throughout the paper. We write $g \in \mathbb{C}_{b}^{k}(\R^{m})$ if the function $g: \R^{m} \to \R$ is bounded and admits continuous bounded partial derivatives up to the order $k$. Recall also the definition of the total variation distance $d_{TV}$ given in (\ref{e:tv}).
\begin{definition}
\begin{enumerate}

\item For every function $g: \R^{m} \to \R$, let $$\Vert g \Vert_{{Lip}}:= \underset{x \neq y}{\sup }\frac{\vert g(x) - g(y) \vert}{\Vert x - y \Vert_{\R^{m}}},$$ where $\Vert \cdot\Vert_{\R^{m}}$ is the usual Euclidian norm on $\R^{m}$. 


\item For a positive integer $k$ and a function $g\in \mathbb{C}_{b}^{k}(\R^{m})$ , we set
$$  \|g^{(k)}\|_{\infty}=\max\limits_{1\leq i_1 \leq \ldots \leq i_k\leq m} \sup\limits_{x\in \R^{m}} \left| \cfrac{\partial^k}{\partial x_{i_1} \ldots \partial x_{i_k}} g(x)\right| .$$
In particular, by specializing this definition to $g^{(2)} = g''$ and $g^{(3)} = g'''$, we obtain
$$  \|g''\|_{\infty}=\max\limits_{1\leq i_1 \leq i_2 \leq m} \sup\limits_{x\in \R^{m}} \left| \cfrac{\partial^2}{\partial x_{i_1}  \partial x_{i_2}} g(x)\right| .$$
$$  \|g'''\|_{\infty}=\max\limits_{1\leq i_1 \leq i_2  \leq i_3\leq m} \sup\limits_{x\in \R^{m}} \left| \cfrac{\partial^3}{\partial x_{i_1} \partial x_{i_2} \partial x_{i_3}} g(x)\right| .$$
\item $\rm{Lip}(1)$ indicates the collection of all real-valued Lipschitz functions, from $\R$ to $\R$, with Lipschitz constant less or equal to one.


\item $\mathscr{C}_{3}$ indicates the collection of all functions $g \in \mathbb{C}_{b}^{3}(\R^{m})$ such that $\Vert g \Vert_{\rm{Lip}} \leq 1$, $\|g''\|_{\infty} \leq 1$ and $\|g'''\|_{\infty} \leq 1$.
\end{enumerate}
\end{definition}

\noindent We now define the different metrics we will use. 

\begin{definition}\label{d:d1}
The metric $d_{\mathscr{H}_{1}}$ between the laws of two $\mathbb{Z}_{+}^{d} \times \R$-- valued random vectors $X$ and $Y$ such that $\E \Vert X \Vert_{\mathbb{Z}_{+}^{d} \times \R}$, $\E \Vert Y \Vert_{\mathbb{Z}_{+}^{d} \times \R} < \infty$, written $d_{\mathscr{H}_{1}}(X,Y)$, is given by
$$d_{\mathscr{H}_{1}}(X,Y) = \underset{h \in \mathscr{H}_{1}}{\rm{sup}}\vert \E(h(X)) - \E(h(Y)) \vert,$$
where $\mathscr{H}_{1}$ indicates the collection of all functions $\psi : \mathbb{Z}_{+}^{d} \times \R \mapsto \R : (j_{1},\ldots,j_{d};x) \mapsto \psi(j_{1},\ldots,j_{d};x)$ such that $\psi$ is bounded by 1 and, for all $j_{1},\ldots,j_{d}$, the mapping $x \mapsto \psi(j_{1},\ldots,j_{d};x)$ is in $\rm{Lip}(1)$.
\end{definition}


\begin{definition}\label{d:d3}
The metric $d_{\mathscr{H}_{3}}$ between the laws of two $\mathbb{Z}_{+}^{d} \times \R^{m}$-- valued random vectors $X$ and $Y$ such that $\E \Vert X \Vert_{\mathbb{Z}_{+}^{d} \times \R^{m}}$, $\E \Vert Y \Vert_{\mathbb{Z}_{+}^{d} \times \R^{m}} < \infty$, written $d_{\mathscr{H}_{3}}(X,Y)$, is given by
$$
d_{\mathscr{H}_{3}}(X,Y) = \underset{h \in \mathscr{H}_{3}}{\rm{sup}}\vert \E(h(X)) - \E(h(Y)) \vert,
$$
where $\mathscr{H}_{3}$ indicates the collection of all functions $\psi : \mathbb{Z}_{+}^{d} \times \R^{m} \mapsto \R : (j_{1},\ldots,j_{d};x_{1},\ldots,x_{m}) \mapsto \psi(j_{1},\ldots,j_{d};x_{1},\ldots,x_{m})$ such that $|\psi|$ is bounded by 1 and for all $j_{1},\ldots,j_{d}$, the mapping $(x_{1},\ldots,x_{m}) \mapsto \psi(j_{1},\ldots,j_{d};x_{1},\ldots,x_{m}) \in \mathscr{C}_{3}$.
\end{definition}

\begin{remark}{\rm The indices $1$ and $3$ label the classes $\mathscr{H}_1$ and $\mathscr{H}_3$, respectively, according to the degree of smoothness of the corresponding test functions. The topology induced by any of the two distances $d_{\mathscr{H}_1}, \, d_{\mathscr{H}_3}$ is strictly stronger than the topology of convergence in distribution. 
}
\end{remark}

\noindent We will sometimes need a useful multidimensional Taylor-type formula on $\mathbb{Z}_+^d$. Given a function $f$ on $\mathbb{Z}_+$, we write $\Delta f(k) = f(k+1)-f(k)$, $k=0,1,...$, and also $\Delta^2 f = \Delta(\Delta f)$. More generally, given a function $f : \mathbb{Z}^d_+ \to \R$, for every $i,j=1,...,d$ we write $\Delta_if(x^{(1)},...,x^{(d)}) = f(x^{(1)},...,x^{(i)}+1,... ,x^{(d)}) - f(x^{(1)},... ,x^{(d)})$, and $\Delta_{ij}^2 = \Delta_i(\Delta_j f)$. Of course, when $d=1$ one has that $\Delta_1 = \Delta$ and $\Delta^2_{11} = \Delta^2$. The proof of the forthcoming statement makes use of the following result, derived in \cite[Proof of Theorem 3.1]{Pec12} (see also \cite{B87}).  For every $f:\mathbb{Z}_{+} \rightarrow \mathbb{R}$, it holds that, for every $k,a\in \mathbb{Z}_{+}$,
\begin{equation}
\label{taylorGio}
f(k) - f(a) = \Delta f(a)(k-a) + R,
\end{equation}
where $R$ is a residual quantity satisfying 
\begin{equation}
\label{residualbound}
\left\lvert R \right\rvert \leq \frac{\Vert \Delta^{2}f\Vert_{\infty}}{2} \left\lvert ( k-a )\left( k-a - 1 \right)  \right\rvert.
\end{equation}
For the rest of the paper, we will use the following notation, which is meant to improve the readability of the proofs. If $ x = \left( x^{(1)},\ldots,x^{(d)}\right) $ is an $d$--dimensional vector, for $k\leq p$ we will denote by $x^{(k,p)}$ the sub-vector composed of the $k$th trough the $p$th component of $x$, i.e. $ x^{(k,p)}= \left( x^{(k)},\ldots,x^{(p)}\right)$. Also, we set by convention $x^{(j,j-1)}  = \emptyset$ for every value of $j$.

\subsection{Complete statement and proof of the Portmanteau inequalities}\label{s:mainstatements}

\noindent We provide below a precise statement of Theorem \ref{t:mainboundintro}, including a discussion of the different cases in terms of dimensions and covariance matrices, each having its own associated metric. The technique of the proof is reminiscent of the computations contained in the classic paper \cite{AGG}. For an explicit description of the constants $K, K_i$, see Remark \ref{r:constants}. 

\medskip

\begin{remark}{\rm
We do not deal with the cases $d=1,\, m=0$ and $d=0, \, m\geq 1$ since they are already covered, respectively, by \cite[Theorem 3.1]{Pec12},  \cite[Theorem 3.1]{PSTU10} and \cite[Theorem 3.3 and Theorem 4.2]{PZ10}. We could have dealt separately with the case where $m\geq 2$ and $C>0$, by using a multidimensional version of Stein's method on the Poisson space, as done in  \cite[Section 3]{PSTU10}: by doing so, we would have been able to consider test functions that are only twice differentiable, as well as bounding constants nicely depending on the operator norm of the matrices $C$ and $C^{-1}$. There is no additional difficulty in implementing this approach (albeit a considerable amount of additional notation should be introduced), and we have refrained to do so merely to keep the length of the paper within bounds. Finally, the results of \cite{PSTU10, PZ10} imply that, in the case $d=0, \, m\geq 1$, one can drop the boundedness assumption for the test functions defining the  distances $d_{\mathscr{H}_1},\,  d_{\mathscr{H}_3}$, as well as the Lipschitz assumption for the functions composing the class $\mathscr{C}_3$. The forthcoming proof will reveal that these boundedness and Lipschitz properties are needed in order to deal with cross terms, that is, expectations involving both elements of ${\bf F}_d$ and ${\bf G}_m$. }
\end{remark}

\begin{theorem}[Portmanteau inequalities: full statement]
\label{mainStatement}
Let $d,m$ be integers such that $d \vee m \geq 1$. Let $H = ({\bf{X}}_{d},{\bf{N}}_{m})$ and $V = ({\bf{F}}_{d},{\bf{G}}_{m})$ be the $(d+m)$--dimensional random elements defined by (\ref{e:H}) and (\ref{e:V}) respectively. Then, the following two statements hold: 
\\~\\
\noindent \textbf{ Case 1}: $d , m \geq 1$. Consider the distance $d_{\mathscr{H}_{1}}$ and $d_{\mathscr{H}_{3}}$, respectively, according as  $m=1$ or $m\geq 2$. For $i=1,3$, there exists a universal positive constant $K_i$ (having the form described in Remark \ref{r:constants})  such that  
\begin{equation}\label{e:goal}
d_{\mathscr{H}_{i}}\left(V,H \right) \leq  K_{i}\left\{ \alpha_1({\bf \lambda}_d, {\bf F}_d)+\alpha_2({\bf F}_d)+\alpha_3({\bf F}_d) + \beta({\bf F}_d,{\bf G}_m) + \gamma_1(C, {\bf G}_m)+\gamma_2({\bf G}_m)\right \}. 
\end{equation}
\\~\\
\noindent \textbf{ Case 2}: $d  \geq 2, \, m=0$. In this case $V = {\bf F}_d$ and $H ={\bf X}_d$, and one has that, for some universal positive constant $K$ (having the form described in Remark \ref{r:constants}),
\begin{equation*}
d_{TV}\left({\bf F}_d ,{\bf X}_d \right) \leq  K\left\{ \alpha_1({\bf \lambda}_d, {\bf F}_d)+\alpha_2({\bf F}_d)+\alpha_3({\bf F}_d) \right \}. 
\end{equation*}

\end{theorem}
Before starting the proof of Theorem \ref{mainStatement}, we prove the following lemma, that is discrete version of Taylor's formula and that will prove quite useful in the upcoming proof of Theorem \ref{mainStatement}.
\begin{lemma}\label{l:multitaylor} Let $f : \mathbb{Z}^d_+ \to \R$. Then, for every $x = (x^{(1)},...,x^{(d)}), \, a = (a^{(1)},...,a^{(d)}) \in \mathbb{Z}^d_+$,
\begin{equation}\label{e:paga}
f(x) = f(a) +\sum_{i=1}^d \Delta_i f(a)(x^{(i)} - a^{(i)}) + R,
\end{equation}
where the residual quantity $R$ satisfies
$$
|R|\leq \frac12 \max_{i,j = 1,...,d} \| \Delta^2_{ij} f\|_\infty \times \left\{ \sum_{i=1}^d |x^{(i)} - a^{(i)}| |x^{(i)} - a^{(i)}-1| + \sum_{1\leq i\neq j \leq d } |x^{(i)} - a^{(i)}| |x^{(j)} - a^{(j)}| \right\}.  
$$
Moreover, one has also the first order estimate:
\begin{equation}\label{e:nini}
\left| f(x) - f(a)\right| \leq \max_{i=1,...,d} \|\Delta_if\|_{\infty} \times \sum_{i=1}^d |x^{(i)} - a^{(i)}|.
\end{equation}

\end{lemma}
\noindent{\textbf{Proof.}
Using (\ref{taylorGio}), one has that
\begin{eqnarray*}
 f(x) - f(a) &=& \sum_{i=1}^d \{ f(a^{(1,i-1)}, x^{(i,d)}) -  f(a^{(1,i)}, x^{(i+1,d)})\}\\
 &=& \sum_{i=1}^d  \Delta_if(a^{(1,i)}, x^{(i+1,d)} ) (x^{(i)} - a^{(i)}) +R_0,
\end{eqnarray*}
where
$$
|R_0|\leq \frac12 \max_{i=1,...,d} \| \Delta^2_{ii} f\|_\infty \times \sum_{i=1}^d |x^{(i)} - a^{(i)}| |x^{(i)} - a^{(i)}-1|.
$$
On the other hand,
$$
\sum_{i=1}^d  \Delta_if(a^{(1,i)}, x^{(i+1,d)} ) (x^{(i)} - a^{(i)})\! =\! \sum_{i=1}^d \Delta_i f(a)(x^{(i)} - a^{(i)}) + \sum_{i=1}^{d-1} \big[ \Delta_if(a^{(1,i)}, x^{(i+1,d)}) -  \Delta_i f(a) \big] (x^{(i)} - a^{(i)}),
$$
and formula (\ref{e:paga}) is immediately obtained from the representation
\begin{eqnarray*}
&&\sum_{i=1}^{d-1} \big[ \Delta_if(a^{(1,i)}, x^{(i+1,d)} )-  \Delta_i f(a) \big] (x^{(i)} - a^{(i)}) \\
&&=\sum_{i=1}^{d-1} (x^{(i)} - a^{(i)}) \sum_{j=i+1}^d \big[ \Delta_i f(a^{(1,j-1)}, x^{(j,d)}) - \Delta_i f(a^{(1,j)}, x^{(j+1,d)})\big],
\end{eqnarray*}
as well as from the elementary inequality
$$
\big| \Delta_i f(a^{(1,j-1)}, x^{(j,d)}) - \Delta_i f(a^{(1,j)}, x^{(j+1,d)})\big| \leq  \| \Delta^2_{ij} f\|_\infty\times | x^{(j)} - a^{(j)}|.
$$
Formula (\ref{e:nini}) follows from 
$$
| f(x) - f(a) | \leq \sum_{i=1}^d | f(a^{(1,i-1)}, x^{(i,d)}) -  f(a^{(1,i)}, x^{(i+1,d)})| \leq \max_{i=1,...,d} \|\Delta_if\|_{\infty}\sum_{i=1}^d |x^{(i)} - a^{(i)}|.
$$
\qed
\\~\\
\noindent{\textbf{Proof of Theorem \ref{mainStatement}}. First of all, we observe that the conclusion of Case 2 follows from the computations leading to the proof of Case 1, by selecting a test function $\psi \in \mathscr{H}_{i}$ uniquely depending on the first $d$ variables. In what follows, $K$ will denote a positive universal constant that may vary from line to line; by a careful bookkeeping of the forthcoming computations, one sees that such a constant $K$ can  be taken to have the form provided in Remark \ref{r:constants}. 

\noindent Now let $\psi \in \mathscr{H}_{i}$. We want to deduce an upper bound for 
\begin{equation*}
\left| \E\left( \psi\left({\bf{F}}_{d},{\bf{G}}_{m}\right)\right) - \E\left( \psi\left({\bf{X}}_{d},{\bf{N}}_{m}\right)\right) \right|. 
\end{equation*}
We can assess such a quantity in the following way:
\begin{eqnarray}\label{e:a2}
&& | \E\left( \psi\left({\bf{F}}_{d},{\bf{G}}_{m}\right)\right) - \E\left( \psi\left({\bf{X}}_{d},{\bf{N}}_{m}\right)\right)| \\
&& \leq  |\E\left( \psi\left({\bf{F}}_{d},{\bf{G}}_{m}\right)\right) - \E\left( \psi\left({\bf{F}}_{d},{\bf{N}}_{m}\right)\right)| +| \E\left( \psi\left({\bf{F}}_{d},{\bf{N}}_{m}\right)\right) - \E\left( \psi\left({\bf{X}}_{d},{\bf{N}}_{m}\right)\right)|.\notag
\end{eqnarray}
The proof will consist of two main steps. In the first one, we will deal with $\E\left( \psi\left({\bf{F}}_{d},{\bf{N}}_{m}\right)\right) - \E\left( \psi\left({\bf{X}}_{d},{\bf{N}}_{m}\right)\right)$ and in the second one with $\E\left( \psi\left({\bf{F}}_{d},{\bf{G}}_{m}\right)\right) - \E\left( \psi\left({\bf{F}}_{d},{\bf{N}}_{m}\right)\right)$. \\~\\
\textit{Step 1: Controlling the term }$\E\left( \psi\left({\bf{F}}_{d},{\bf{N}}_{m}\right)\right) - \E\left( \psi\left({\bf{X}}_{d},{\bf{N}}_{m}\right)\right)$. 
Such a term can be decomposed in the following way:
\begin{equation*}
\E\left( \psi\left({\bf{F}}_{d},{\bf{N}}_{m}\right)\right) - \E\left( \psi\left({\bf{X}}_{d},{\bf{N}}_{m}\right)\right) = \sum_{k=1}^{d}\E\left(\psi\left({\bf{X}}^{(1,k-1)},{\bf{F}}^{(k,d)},{\bf{N}}_m\right) - \psi\left({\bf{X}}^{(1,k)},{\bf{F}}^{(k+1,d)},{\bf{N}}_m\right)\right).
\end{equation*}
We will now study separately each term appearing in the sum. In what follows, we write $\mathscr{L}_U$ to indicate the probability measure given by the law of a given random element $U$; integrals with respect to $\mathscr{L}_U$ are implicitly realized over the set where $U$ takes values. For any fixed $1\leq k \leq d$, by exploiting independence, we have
\begin{eqnarray*}
&& \E\left[\psi\left({\bf{X}}^{(1,k-1)},{\bf{F}}^{(k,d)},{\bf{N}}_m\right) - \E\left( \psi\left({\bf{X}}^{(1,k)},{\bf{F}}^{(k+1,d)},{\bf{N}}_m\right)\right) \right] = \\
&&\int \mathscr{L}_{{\bf{F}}^{(k,d)}} (dx^{(k,d)})\\ 
&& \quad\quad\quad\quad\E\left\lbrace \psi\left({\bf{X}}^{(1,k-1)},x^{(k,d)},{\bf N}_m \right) \!\!-\!\! \int \mathscr{L}_{X^{(k)}}(da)\psi\left({\bf{X}}^{(1,k-1)},a,x^{(k+1,d)},{\bf N}_m \right)  \right\rbrace .
\end{eqnarray*}
For a fixed $\left(z^{(1,k-1)},x^{(k+1,d)},y \right) \in \mathbb{Z}_{+}^{d-1} \times \R^{m}$, we denote by $$x^{(k)}\mapsto f_{k}\left( z^{(1,k-1)},x^{(k)},x^{(k+1,d)},y \right) := f_k(x^{(k)})$$ the unique solution to the {\it Chen-Stein equation}
\begin{eqnarray*}
\widetilde{\psi}(x^{(k)}) - \E(\widetilde{\psi}(X^{(k)})  )  &=& \lambda_{k}f (x^{(k)}+1 ) - x^{(k)}f(x^{(k)}), \quad x^{(k)} = 0,1,..., 
\end{eqnarray*}
satisfying the boundary condition $\Delta^2f(0) = 0$, where $\widetilde{\psi}(x^{(k)}) :=\psi\left( z^{(1,k-1)},x^{(k)},x^{(k+1,d)},y\right)$. We recall (see e.g. \cite{ErSurvey}) that $f_k$ is given by $f_k(0) = f_k(1) - \Delta f_k(2)$ and, for $x =1,2,...$,
\begin{eqnarray}\label{e:pimpa}
f_k(x )& =& \frac{(x - 1)!}{\lambda^x_k} \sum_{w=0}^{x-1}\left[ \frac{\lambda_k^w}{w!}\left(\widetilde{\psi}(w) - \E[\widetilde{\psi}(X^{(k)})]\right)  \right] \\
&=& - \frac{(x - 1)!}{\lambda^x_k} \sum_{w=x}^{\infty}\left[ \frac{\lambda_k^w}{w!}\left(\widetilde{\psi}(w) - \E[\widetilde{\psi}(X^{(k)})]\right)   \right].\notag
\end{eqnarray}
Using the fact that $|\widetilde{\psi}|\leq 1$ together with \cite[Theorem 2.3]{ErSurvey} and \cite[Theorem 1.3]{daly}, we deduce that $|f_k| \leq 3$, $|\Delta f_k| \leq 2 (1 - e^{-\lambda_k})/\lambda_k$ and $|\Delta^2 f_k| \leq 4 (1 - e^{-\lambda_k})/\lambda^2_k$.\footnote{The upper bound on $|f_k|$ can be reduced to 2 if one selects a solution of the Chen-Stein equation such that $f(0) = 0$} Exploiting once again independence, we can now write:
\begin{eqnarray}
\label{mainestimatestep1}
&& \E\left(\psi\left({\bf{X}}^{(1,k-1)},{\bf{F}}^{(k,d)},{\bf{N}}_m\right) - \E\left( \psi\left({\bf{X}}^{(1,k)},{\bf{F}}^{(k+1,d)},{\bf{N}}_m\right)\right) \right) =   \\
&& \int \mathscr{L}_{{\bf{F}}^{(k,d)}}(dx^{(k,d)} )\notag \\
&&\quad\quad\quad\E \left\lbrace \lambda_{k}f_{k}\left( {\bf{X}}^{(1,k-1)},x^{(k)}+1,x^{(k+1,d)},{\bf{N}}_m \right) - x^{(k)}f_{k}\left(  {\bf{X}}^{(1,k-1)},x^{(k)},x^{(k+1,d)},{\bf{N}}_m \right)  \right\rbrace \nonumber \\
&=&\!\!\!\! \E\left(\lambda_{k}f_{k}\left( {\bf{X}}^{(1,k-1)},F^{(k)}+1,{\bf{F}}^{(k+1,d)},{\bf{N}}_m\right) - F^{(k)}f_{k}\left(  {\bf{X}}^{(1,k-1)},F^{(k)},{\bf{F}}^{(k+1,d)},{\bf{N}}_m\right) \right)\nonumber  \\
&=&\!\!\!\! \E\left(\lambda_{k}\Delta f_{k}\left( {\bf{X}}^{(1,k-1)},{\bf{F}}^{(k,d)},{\bf{N}}_m \right) - \delta\left(-DL^{-1}F^{(k)} \right)f_{k}\left(  {\bf{X}}^{(1,k-1)},{\bf{F}}^{(k,d)},{\bf{N}}_m\right) \right)\nonumber  \\
&=&\!\!\!\! \E\left(\!\lambda_{k}\Delta f_{k}\left( {\bf{X}}^{(1,k-1)},{\bf{F}}^{(k,d)},{\bf{N}}_m \right) \!- \! \left\langle \!Df_{k}\!\left(  {\bf{X}}^{(1,k-1)},{\bf{F}}^{(k,d)},{\bf{N}}_m\right) , -DL^{-1}F^{(k)} \right\rangle_{L^{2}(\mu)}   \right).\notag
\end{eqnarray}
Note that (since $H$ is assumed to be independent of $\eta$) in the previous expressions the Malliavin operators act on random variables only through their dependence on the components of ${\bf F}_d$. We now need to explicitly calculate $Df_{k}\left(  {\bf{X}}^{(1,k-1)},{\bf{F}}^{(k,d)},{\bf{N}}_m\right)$, and (by virtue of (\ref{e:diffop})), one has
\begin{eqnarray}
&& D_{z}f_{k}\left(  {\bf{X}}^{(1,k-1)},{\bf{F}}^{(k,d)},{\bf{N}}_m\right) \nonumber \\
\label{partiTwo}
&& \quad =  f_{k}\left({\bf{X}}^{(1,k-1)},{\bf{F}}_{z}^{(k,d)},{\bf{N}}_m\right) - f_{k}\left(  {\bf{X}}^{(1,k-1)},{\bf{F}}^{(k,d)},{\bf{N}}_m\right).
\end{eqnarray}
In order to deal with this quantity, one should first observe that, for every $k$, the mapping $f_k(\cdot, {\bf N}_m) : \mathbb{Z}^{d} \to \R$ given by $$x \mapsto f_k \left( x,{\bf{N}}_m\right)$$ takes values in $[-3,3]$, and therefore $\| \Delta_{i} f_k(\cdot, {\bf N}_m)\|_\infty \leq 6 $ and $\| \Delta^2_{i,j} f_k(\cdot, {\bf N}_m)\|_\infty \leq 12 $, for every $i,j=1,...,d$. One can now use Lemma \ref{l:multitaylor} to deduce that
$$
D_{z}f_{k}\left(  {\bf{X}}^{(1,k-1)},{\bf{F}}^{(k,d)},{\bf{N}}_m\right) =\sum_{i=k}^d \Delta_if_k\left(  {\bf{X}}^{(1,k-1)},{\bf{F}}^{(k,d)},{\bf{N}}_m\right)D_{z}F^{(i)} +R^{(k)}_z,
$$
where 
$$
|R_z^{(k)}|\leq 6 \times \left\{ \sum_{i=k}^d |D_{z}F^{(i)} | |D_{z}F^{(i)}-1| + \sum_{k\leq i\neq j \leq d } |D_{z}F^{(i)}| |D_{z}F^{(j)}| \right\}.
$$
Using the fact that (by definition) $$\Delta_k f_k\left(  {\bf{X}}^{(1,k-1)},{\bf{F}}^{(k,d)},{\bf{N}}_m\right) = \Delta f_k\left(  {\bf{X}}^{(1,k-1)},{\bf{F}}^{(k,d)},{\bf{N}}_m\right),$$ gathering the previous estimates together, and applying them to (\ref{mainestimatestep1}) finally gives:
$$
\left\lvert \E\left( \psi\left({\bf{F}}_d,{\bf{N}}_m\right)\right) - \E\left( \psi\left({\bf{X}}_d,{\bf{N}}_m\right)\right) \right\rvert \leq K\left\lbrace \alpha_{1}(\lambda_d,{\bf{F}}_d) + \alpha_{2}({\bf{F}}_d) + \alpha_{3}({\bf{F}}_d) \right\rbrace .
$$

\noindent \textit{Step 2: Controlling the term }$\E\left( \psi\left({\bf{F}}_d,{\bf{G}}_m\right)\right) - \E\left( \psi\left({\bf{F}}_d,{\bf{N}}_m\right)\right)$. 
This part is slightly more delicate, since one has to take into account the dependence between ${\bf F}_d$ and ${\bf G}_m$. We have to consider two cases, namely $m=1$ and $m\geq 2$. Note that, in the second case, it is not necessary to assume that the matrix $C$ is positive definite.
~\\~\\
\textit{$(m=1)$ }
In this case ${\bf G}_m$ and ${\bf N}_m$ are two real-valued random variables $G\in {\rm dom}D$ and $N\sim \mathscr{N}(0,1)$. We will only consider the case $C=1$, and one can recover the general statement by elementary considerations. For every $x\in \mathbb{Z}^d_+$ and $y\in \R$, we write 
$$
f_\psi(x,y) = e^{y^2/2} \int_{-\infty}^y \{\psi(x,a) - \E[\psi(x,N)]\} e^{-a^2/2}da.
$$
It is well-known (see e.g. \cite[Chapter 3]{np-book}) that $f_\psi$ satisfies the (parametrized) Stein equation
$$
\partial_yf_\psi(x,y) - yf_\psi(x,y) = \psi(x,y) - \E[\psi(x,N)], \quad y\in \R, \, x\in \mathbb{Z}_+^d,
$$
where we have used the symbol $\partial_y$ to indicate a partial derivative with respect to $y$. Moreover, thanks to the assumptions on $\psi$, one can prove that the following relations are in order for every $x$: $\|f_\psi(x,\cdot)\|_\infty \leq \sqrt{2\pi}$, $\|\partial_y f_\psi(x,\cdot)\|\leq 1$, and $\|\partial^2_{yy} f_\psi(x,\cdot)\|_\infty \leq 2$ (note that the partial derivatives $\partial^2_{yy} f_\psi(x,\cdot)$ are only defined up to a subset of $\R$ of measure 0). It follows that
\begin{eqnarray}
&& \E\left( \psi \left( {\bf F}_d,G \right)\right) - \E\left( \psi\left({\bf F}_d,N\right)\right) = \E[\partial_yf_\psi({\bf F}_d,G) - Gf_\psi({\bf F}_d,G)]\notag \\
&& = E[\partial_yf_\psi({\bf F}_d,G)] - \E[\langle -DL^{-1}G , Df_\psi({\bf F}_d,G)\rangle_{L^2(\mu)} ].\label{e:ee}
\end{eqnarray}
Clearly, $D_z f_\psi({\bf F}_d,G) = A_z+B_z$, where
$$
A_z := f_\psi(({\bf F}_d)_z, G_z) - f_\psi({\bf F}_d, G_z), \quad B_z := f_\psi({\bf F}_d, G_z) - f_\psi({\bf F}_d, G)
$$
Using a Taylor expansion as in \cite[Proof of Theorem 3.1]{PSTU10}, one sees that
$$
B_z = \partial_y f_\psi({\bf F}_d, G)D_{z}G +R_z, 
$$
where $|R_z| \leq (D_zG)^2$. Now observe that the mapping $f_\psi( \cdot, G_z) : \mathbb{Z}^d \to \R : x \mapsto f_\psi( x, G_z)$ is bounded by $\sqrt{2\pi}$, in such a way that $\| \Delta_i f_\psi( \cdot, G_z)\|_\infty \leq 2\sqrt{2\pi}$, for every $i=1,...,d$. We can therefore use formula (\ref{e:nini}) to infer that
$$
| A_z |  \leq 2\sqrt{2\pi} \sum_{i=1}^d | D_{z}F^{(i)}|.
$$ 
Plugging these relations into (\ref{e:ee}) yields that 
$$
\left|  \E\left( \psi \left( {\bf F}_d,G \right)\right) - \E\left( \psi\left({\bf F}_d,N\right)\right) \right| \leq K\{\beta({\bf F}_d,G) + \gamma_1(C, G)+\gamma_2(G)\}.
$$
~\\
\textit{$(m\geq 2)$ } We use an interpolation technique analogous to the one appearing in \cite[Proof of Theorem 4.2]{PZ10}. For every $t\in [0,1]$, we define
$$
\Phi(t) := \E\{ \psi({\bf F}_d, \sqrt{1-t}{\bf G}_m +\sqrt{t} {\bf N}_m)\},
$$
in such a way that $\left| \E\{ \psi({\bf F}_d, {\bf G}_m)\} - \E\{ \psi({\bf F}_d, {\bf N}_m)\}\right| \leq \int_0^1 | \Phi'(t)|dt$.  Taking the derivative with respect to $t$ and then integrating by parts shows that
$$
\Phi'(t) =A_t - B_t,
$$
where, with obvious notation,
$$
A_t = \frac12 \sum_{i,j=1}^d C(i,j) \E\left[ \partial^2_{y_iy_j} \psi({\bf F}_d, \sqrt{1-t}{\bf G}_m +\sqrt{t} {\bf N}_m)\right],
$$
and
\begin{eqnarray*}
&&B_t = \frac{1}{2\sqrt{1-t}}\sum_{j=1}^m \E\left[ \langle -DL^{-1}G^{(j)}, D\partial_{y_j} \psi({\bf F}_d, \sqrt{1-t}{\bf G}_m +\sqrt{t} {\bf N}_m) \rangle_{L^2(\mu)} \right]\\
&& = \frac{1}{2\sqrt{1-t}}\sum_{j=1}^m\left\{ \E\left[ \langle -DL^{-1}G^{(j)}, b^{1,j} \rangle_{L^2(\mu)} \right] + \E\left[ \langle -DL^{-1}G^{(j)}, b^{2,j} \rangle_{L^2(\mu)} \right]\right\}\\
&&:= B^{(1)}_t +B_t^{(2)},
\end{eqnarray*}
where the random functions $b^{1,j}$ and $b^{2,j}$ are given by 
$$
z \mapsto b^{1,j}_z:= \partial_{y_j}\psi(({\bf F}_d)_z, \sqrt{1-t}({\bf G}_m)_z +\sqrt{t} {\bf N}_m) - \partial_{y_j}\psi({\bf F}_d, \sqrt{1-t}({\bf G}_m)_z +\sqrt{t} {\bf N}_m), 
$$
and 
$$
z \mapsto b^{2,j}_z:=\partial_{y_j} \psi({\bf F}_d, \sqrt{1-t}({\bf G}_m) _z+\sqrt{t} {\bf N}_m) - \partial_{y_j}\psi({\bf F}_d, \sqrt{1-t}{\bf G}_m +\sqrt{t} {\bf N}_m).
$$
Reasoning exactly as in the proof of \cite[Theorem 4.1]{PZ10}, one proves that 
$$
\sup_{t\in [0,1]} \left| A_t - B^{(2)}_t\right| \leq \frac{1}{4} \{ \gamma_1(C, {\bf G}_m)+\gamma_2({\bf G}_m)\}.
$$
To conclude, we apply again Lemma \ref{l:multitaylor}. Start by observing that, since 
$|\partial_{y_j} \psi | \leq 1 $ by assumption, one has that, for every $i=1,...,d$, $\| \Delta_i \partial_{y_j} \psi(\cdot, \sqrt{1-t}({\bf G}_m) _z+\sqrt{t} {\bf N}_m)\|_\infty\leq 2$. We can now use (\ref{e:nini}) to infer that
$$
| b^{1,j}_z | \leq 2 \sum_{i=1}^d | D_zF^{(i)}|.
$$ 
These estimates yield eventually that
$$
\int_0^1 |B_t^{(1)}|dt \leq \int_0^1 \frac{2}{\sqrt{1-t}}dt \times \beta({\bf F}_d, {\bf G}_m) = 4  \beta({\bf F}_d, {\bf G}_m),
$$
and the desired conclusion follows at once.

\qed

\subsection{Proof of Theorem \ref{stablePoisson}}\label{ss:poissonproof}

The first part of the statement is the same as Case 2 of Theorem \ref{mainStatement}. In order to deduce the conclusion about stable convergence, one should fix an integer $l\geq 1$, as well as pairwise disjoint sets $A_1,...,A_l\in \mathscr{Z}_\mu$, and then build an ancillary $(d+l)$-dimensional vector 
$$
{\bf F}'_{d+l,n} := ({\bf F}_{d,n},\,  \eta(A_1),...,\eta(A_l) ).
$$
Applying again Case 2 of Theorem \ref{mainStatement}, one proves immediately that conditions (\ref{e:aaa})--(\ref{e:cat}) imply that ${\bf F}'_{d+l,n}$ converges in distribution to $\left({\bf X}_d, \hat\eta(A_1),...,\hat\eta(A_l)\right)$. Since ${\bf X}_d$ is independent of $\eta$ by definition, we deduce that, for every $(\gamma_1,...,\gamma_d)\in \mathbb{R}^d$, every collection $A_1,...,A_l\in \mathscr{Z}_\mu$ of disjoint sets and every random variable $Z = \varphi(\eta(A_1),...,\eta(A_l))$ with $\varphi$ bounded,
$$
\underset{n \rightarrow \infty}{\rm lim}\E\left[e^{i\langle {\bf F}_{d,n}, \gamma\rangle_{\R^d}} Z\right] = \E[Z]\times \E\left[e^{i\langle {\bf X}_d, \gamma\rangle_{\R^d}} \right].
$$
An application of Point 4 of Lemma \ref{l:stable} yields the desired conclusion.
\qed\\

\subsection{Proof of Theorem \ref{perturbedStablePoisson}}\label{ss:stpoissmwip}

\noindent {\it Step 1: convergence in distribution}. We start by proving that ${\bf F}_{d,n}$ converges in distribution to ${\bf X}_d$. Our plan is to apply Case 2 of Theorem \ref{mainStatement}. Exploiting the fact that each $\{B_n^{(i)}\}$ is a smooth vanishing perturbation, and reasoning exactly as in the first part of the proof of \cite[Theorem 4.12]{Pec12}, one sees that it is enough to prove that Conditions 1 and 2 imply that the five sums appearing in the definitions of $\alpha_1(\cdot), \alpha_2(\cdot), \alpha_3(\cdot)$ (see (\ref{e:alpha1})--(\ref{e:alpha3})) all converge to zero, whenever one chooses the vector $\big( I_{q_1} (f_n^{(1)}),...,I_{q_d} (f_n^{(d)}) \big)$ as their argument. Again from the proof of \cite[Theorem 4.12]{Pec12}, we know that the assumptions in the statement imply that, for every $i=1,...,d$
$$
\underset{n \rightarrow \infty}{\rm lim}\left\{  \E\left[\left| \lambda_i - q_i^{-1}\|DI_{q_i}(f_n^{(i)})\|^2_{L^2(\mu)} \right|\right] + \E\left[ \int_Z (D_zI_{q_i}(f_n^{(i)}))^2(D_zI_{q_i}(f_n^{(i)})-1)^2 \mu(dz) \right]\right\} =0.
$$
Using the fact that the sequence
$$n\mapsto \E\left[ \int_Z \left(D_z I_{q_i}(f_n^{(i)}) \right)^2 \mu(dz)\right]= q_i^2 \E\left[ \int_Z \left(D_z L^{-1}I_{q_i}(f_n^{(i)}) \right)^2 \mu(dz)\right] = q_i\E[(F_n^{(i)})^2] $$
is bounded for every $i$, and by a standard application of the Cauchy-Schwarz inequality, we see that it is enough to prove that, for every $i\neq j$,
\begin{equation}\label{e:lim}
\underset{n \rightarrow \infty}{\rm lim}\left\{  \E \left[ \int_Z (D_zI_{q_i}(f_n^{(i)}))^2(D_zI_{q_j}(f_n^{(j)}))^2 \mu(dz) \right] + \E\left[\langle DI_{q_i}(f_n^{(i)}),DI_{q_j}(f_n^{(j)})\rangle^2_{L^2(\mu)} \right] \right\}  =0.
\end{equation}
Using the computations contained in \cite[p. 464]{PSTU10}, one sees that, for every $i$,
\begin{equation}\label{e:dev}
\{D_zI_{q_i}(f^{(i)}_n)\}^2 =q_i^2 \sum_{p=0}^{2q_i-2} I_p(G_p^{q_i-1}f_n^{(i)}(z,\cdot)),
\end{equation}
where 
\begin{equation}\label{e:der2}
G_{p}^{q_i-1}f_n^{(i)}(z,\cdot)(z_1,...,z_p) = \sum_{r=0}^{q_i-1}\sum_{l=0}^r {\bf 1}_{\{2q_i-2-r-l=p\}} r!\binom{q_i-1}{r}^2\binom{r}{l} \widetilde{f_n^{(i)}(z,\cdot)\star_r^lf_n^{(i)}(z,\cdot)}(z_1,...,z_p),
\end{equation}
where the tilde indicates a symmetrization with respect to the variables represented by a dot (in such a way that the symmetrization does not involve the variable $z$), and the stochastic integrals are set equal to zero on the exceptional set composed of those $z$ such that $f_n^{(i)}(z,\cdot)\star_r^lf_n^{(i)}(z,\cdot)$ is not an element of $L^2(\mu^{2q_i-2-r-l})$ for some $r,l$. We can assume without loss of generality that $q_i\leq q_j$. Applying the isometric properties of multiple integrals using the Fubini theorem and integrating over $Z$, we see that the first summand in (\ref{e:lim}) is a linear combination of objects of the type
$$
C_n = C_n(l,r,s,t,p) :=\int_Z \langle  \widetilde{f_n^{(i)}(z,\cdot)\star_r^lf_n^{(i)}(z,\cdot)} ,  \widetilde{f_n^{(j)}(z,\cdot)\star_t^s f_n^{(j)}(z,\cdot)} \rangle_{L^2(\mu^p)} \mu(dz),
$$
where the indices satisfy the following constraints: (i) $p=0,... ,2q_i -1$, (ii) $r=0,...,q_i-1$, (iii) $l=0,...,r$, (iv)$t=0,...,q_j-1$, (v) $=0,...,t$, and (vi) $2q_i-2-r-l = 2q_j-2-t-s = p$. In the case where $q_i = q_j = r=t$, $l=s=0$ (and therefore $p = q_i-1$), one has that $C_n = \Big \langle ( f_{n}^{(i)} )^2,(f_{n}^{(j)})^2\Big \rangle_{L^2(\mu^{q_i})} \to 0$.  In all other cases, one can prove that
$$
|C_n| \leq \sqrt{\int_Z  \|  {f_n^{(i)}(z,\cdot)\star_r^lf_n^{(i)}(z,\cdot)}  \|_{L^2(\mu^p)} \mu(dz)} \times
 \sqrt{\int_Z  \|  {f_n^{(j)}(z,\cdot)\star_s^t f_n^{(j)}(z,\cdot)}  \|_{L^2(\mu^p)} \mu(dz)}\rightarrow 0
$$ 
(where the first inequality follows from the Cauchy-Schwarz inequality) by directly applying the computations contained in \cite[p. 467]{PSTU10}, as well as the fact that (by assumption) $\underset{n}{\rm sup\ }  \|f_n^{(i)}\|_{L^4(\mu^{q_i} )}<\infty$ for every $i$. We now focus on the second summand in (\ref{e:lim}). We can use directly \cite[Proposition 5.5]{PZ10} to deduce that, whenever $q_i = q_j$ the quantity $ \E\left[\langle DI_{q_i}(f_n^{(i)}),DI_{q_j}(f_n^{(j)})\rangle^2_{L^2(\mu)} \right]$ is equal to a finite linear combination of the squared inner product $\langle f_n^{(i)}, f_n^{(j)} \rangle_{L^2(\mu^{q_i})}^2$, as well as of products of norms of the type
\begin{equation}\label{e:none}
\| f_n^{(i)} \star^{q_i - t}_{q_i - s(t,k)} f_n^{(i)}\|_{L^2(\mu^{t+s(t,k)})} \times \| f_n^{(j)} \star^{q_j - t}_{q_j - s(t,k)} f_n^{(j)}\|_{L^2(\mu^{t+s(t,k)})},
\end{equation}
where $s(t,k) = 2q_i -k - t$ and the indices satisfy the constraints: $k=1,...,2q_i-2$, $t=1,...,q_i$ and $1\leq s(t,k) \leq t$. On the other hand, when $q_i\neq q_j$ the same Proposition 5.5 in \cite{PZ10} implies that $ \E\left[\langle DI_{q_i}(f_n^{(i)}),DI_{q_j}(f_n^{(j)})\rangle^2_{L^2(\mu)} \right]$ is a finite linear combination of products of norms of the type (\ref{e:none}), where $s(t,k) = q_i+q_j -k - t$ and the indices satisfy the constraints: $k=|q_i-q_j|,...,q_i+q_j-2$, $t=1,...,q_i\wedge q_j$ and $1\leq s(t,k) \leq t$. In both cases, the involved products of norms converge to zero whenever Condition 1 in the statement is satisfied, and Case 2 of Theorem \ref{mainStatement} implies that ${\bf F}_{d,n}$ converges in distribution to ${\bf X}_d$ in the sense of total variation.

\medskip

\noindent {\it Step 2: stable convergence.} We apply the second part of Theorem \ref{stablePoisson}. In view of the previous computations, and by reasoning again as at the beginning of the previous step, it is enough to show that, for every $A\in \mathscr{Z}_\mu$ and every $i=1,...,d$,
$$
\E\left[ \left( \int_{A} D_{z}I_{q_i}(f_n^{(i)}) \mu(dz) \right)^2\right] \to 0.
$$
This follows immediately from the relation
$$
\int_{A} D_{z}I_{q_i}(f_n^{(i)}) \mu(dz) = q_i I_{q_i-1} (f_n^{(i)} \star_1^1 g)
$$
where $g(z) = {\bf 1}_A(z)$, as well as $\| f_n^{(i)} \star_1^1 g\|^2_{L^2(\mu^{q_i-1})} = \langle f_n^{(i)}\star_{q_i-1}^{q_i-1} f_n^{(i)}, g\star_0^0g \rangle_{L^2(\mu^{2})}$ (which follows from a Fubini argument).

\qed

\subsection{Proofs of Theorem \ref{perturbedStableGaussian} and Proposition \ref{p:asint} }\label{ss:mwicltp}

\noindent {\it Proof of Theorem \ref{perturbedStableGaussian}.} In view of \cite[Theorem 5.8]{PZ10} we only need to prove stable convergence. To do this we fix an integer $d\geq 1$, as well as disjoint sets $A_1,...,A_d\in \mathscr{Z}_\mu$. Using as before Point 4 of Lemma \ref{l:stable}, the desired conclusion is achieved if we show that the $(d+m)$-dimensional vectors $({\bf F}_d, {\bf G}_{m,n})$, $n\geq 1$ where ${\bf F}_d= (\eta(A_1),...,\eta(A_d))$, converge in distribution to $( {\bf F}_d, {\bf N}_m)$ (recall that ${\bf N}_m$ is independent of $\eta$ by definition). Define ${\bf \lambda}_d = (\mu(A_1),...,\mu(A_d))$.  One has that $\alpha_1({\bf \lambda}_d, {\bf F}_d)+\alpha_2({\bf F}_d)+\alpha_3({\bf F}_d) = 0$, and also that, under the assumptions in the statement, $\gamma_1(C, {\bf G}_{m,n})+\gamma_2({\bf G}_{m,n}) \to 0$ (as a consequence of \cite[Theorem 5.8]{PZ10}). To conclude, we have to show that $ \beta({\bf F}_d,{\bf G}_{n,m}) \to 0$. This follows immediately from Proposition \ref{prophelperalpha3}, since the computations contained in \cite[Proof of Theorem 5.1]{PSTU10} imply that, under the assumptions in the statement,
\begin{equation}\label{e:jog}
\E\left[  \int_{A} [ D_{z}I_{q_i}(g_n^{(j)})]^4 \mu(dz) \right] \to 0, \quad \forall j=1,...,m.
\end{equation}
\qed

\noindent {\it Proof of Proposition \ref{p:asint}}. In view of Proposition \ref{prophelperalpha3}, the conclusion is an immediate consequence of relation (\ref{e:jog}).\qed

\subsection{Proofs of Theorem \ref{t:indu} and Proposition \ref{t:indue}}

\noindent {\it Proof of Theorem \ref{t:indu}}. For every $n$, let $X_{\lambda_n}$ be a one-dimensional Poisson random variable of parameter $\lambda_n$, and recall (see e.g. \cite{AL}) that $d_{TV}(X_\lambda, X_{\lambda_n})\leq |\lambda - \lambda_n|$. The distance $d_\star$ in the statement can be chosen to be $d_{\mathscr{H}_1}$ (see Definition \ref{d:d1}). An application of the triangular inequality and of the independence between $X_\lambda$ and $N$ yield that
$$
d_\star(V_n, H) \leq d_{TV}(X_\lambda, X_{\lambda_n}) + d_\star(V_n, (X_{\lambda_n}, N)).
$$
The conclusion follows from Theorem \ref{mainStatement}, since one has that:
\begin{itemize}

\item[--] according to \cite[Proof of Proposition 4.1]{SchTh2012}, $|\lambda - \lambda_n| + \alpha_1(\lambda_n, F_n) +\alpha_2(F_n) \leq A_n$;

\item[--] according to \cite{LacPec_a}, $\gamma_1(1,\tilde{G}_n) +\gamma_2(\tilde{G}_n) \leq C_kB(G_n, \sigma_n)$;

\item[--] by virtue of the H\"older inequality, and of the fact that $F_n$ takes values in $\mathbb{Z}_+$, 
$$
\beta(F_n, \tilde{G}_n) \leq \! \E\left[ \int_Z \left|D_z F_n \right|^2 \!\! \mu_n(dz)\right]^{\frac{3}{4}} \!\!\times \E\left[ \int_Z \left|D_z L^{-1}\tilde {G}_n \right|^{4} \!\! \mu_n(dz)\right]^{\frac{1}{4}}\leq R \times B(G_n; \sigma_n)^{1/2},
$$
for some constant $R$ independent of $n$, where we have used the fact that, since $F_n$ and $G_n$ both live in a finite sum of Wiener chaoses (see Proposition \ref{p:zap}), then (a) the mapping $n\mapsto \E\left[ \int_Z \left|D_z F_n \right|^2 \!\! \mu_n(dz)\right]$ is bounded, and (b) $$\E\left[ \int_Z \left|D_z L^{-1}\tilde {G}_n \right|^{4} \!\! \mu_n(dz)\right] \leq \E\left[ \int_Z \left|D_z \tilde {G}_n \right|^{4} \!\! \mu_n(dz)\right]\leq C_k B(G_n; \sigma_n)^{2}$$ for every $n$. 
\end{itemize}
\qed 

\medskip

\noindent{\it Proof of Proposition \ref{t:indue}}. In view of the standard theory of Hoeffding decompositions (see e.g. \cite{Vitale}), for every $n\geq 1$, both $F^0_n$ and $\tilde{G}^0_n$ have the form of a $U$-statistic of the type 
\begin{equation}\label{e:unno}
U_n = \E[U_n] + \sum_{l=1}^m \sum_{\{i_1,...,j_l\}\subset [n]} U_{n,l}(Y_{j_1},...,Y_{j_l}),
\end{equation}
where $[n]= \{1,...,n\}$, $m$ is the order of the $U$-statistic (so, $m=k$ or $m=k'$, according as $U_n = \tilde{G}_n$ or $U_n = F_n$), and every kernel $U_{n,l}$ is a symmetric function in $l$ variables satisfying the Hoeffding-type degeneracy condition: $\E[U_{n,l}(Y_1,...,Y_l)| Y_1,...,Y_{l-1}] = 0$. The mean and variance of $F_n$ and $\tilde{G}_n$ are both converging, and this implies that, since the mapping
$$
n\mapsto E[U^2_n] = \E[U_n]^2 + \sum_{l=1}^n \binom{n}{l} \E[U_{n,l}(Y_{1},...,Y_{l})^2]
$$
converges to a finite limit, then the sequences $ n\mapsto \binom{n}{l} \E[U_{n,l}(Y_{1},...,Y_{l})^2]$, $l=1,...,m$, are necessarily bounded. Now, it is easily seen that
$U_n$ is the de-poissonized version of the poissonized $U$-statistic $U'_n$ obtained by replacing $[n]$ with $[N(n)] = \{1,...,N(n)\}$ in the second sum on the RHS of (\ref{e:unno}). The desired conclusion follows from the forthcoming Lemma \ref{l:lalla}, whose proof uses computations from \cite{DyMa}. \qed

\begin{lemma}[De-poissonization Lemma]\label{l:lalla} Let the above notation and assumptions prevail. Then, as $n\to \infty$,
$$
\E[(U_n - U'_n)^2]\rightarrow 0.
$$
\end{lemma}
\noindent{\it Proof.}  Conditioning on $N(n)$, and using standard results on the moments of Poisson random variables, yields (as $n\to\infty$)
$$
E[U'^2_n] = \E[U_n]^2 + \sum_{l=1}^m \E\left[ \binom{N(n)}{l}\right] \E[U_{n,l}(Y_{1},...,Y_{l})^2] \rightarrow c : = \underset{n \rightarrow \infty}{\rm lim}E[U^2_n].
$$
Conditioning again on $N(n)$, we infer that 
$$
E[U_n U'_n] = \sum_{l=1}^m \binom{n}{l} \E[U_{n,l}(Y_{1},...,Y_{l})^2]  b_{n,l},
$$
where $b_{n,l} = \sum_{p=0}^\infty e^{-p}\frac{n^p}{p!} \binom{n\wedge p}{l} \binom{n}{l}^{-1}$. To conclude, it remains to apply the computations contained in \cite[p. 745]{DyMa}, which imply that $b_{n,l} \to 1$ for every $l$.
\qed

\section{Random graphs: proof of Theorem \ref{t:mainrg}}\label{s:fullrg}

\noindent The distance $d_\star$ appearing in the statement is the distance $d_{\mathscr{H}_1}$ introduced in Definition \ref{d:d1}. First of all we observe that, for every $a=0,1,...,d$, the random variable $G'_n(\Gamma_a)$ has the form of a $U$-statistic, that is:
\[
G'_n(\Gamma_a) = \sum_{(x_1,...,x_{k_a}) \in \eta_{n,\neq}^{k_a}} h_{\Gamma_a,t_n}(x_1,...,x_{k_a}),
\]
where: (i) $k_a = k$ for $a=1,...,d$, (ii) the notation indicates that the sum runs over all ordered vectors $(x_1,...,x_{k_a})$ such that each $x_l$ is in the support of $\eta_n$ and $x_l\neq x_{l'}$ for $l\neq l'$, and (iv) the quantity $h_{\Gamma_a,t_n}(x_1,...,x_{k_a})$ is equal to $1/k_a!$ if the restriction of $G'(Y;t_n)$ to $\{x_1,...,x_{k_a}\}$ is isomorphic to $\Gamma_a$ and equals 0 otherwise. It is clear that, for every $a$, the mapping $h_{\Gamma_a,t_n}: (\R^m)^{k_a} \to \R $ is symmetric and {stationary}, in the sense that it only depends on the norms $\|x_l - x_m\|_{\R^m}$, $l\neq m$. We can now apply Proposition \ref{p:zap} to deduce that $G'_n(\Gamma_a) $ admits the following chaotic decomposition
\begin{equation}\label{e:proofc}
G'_n(\Gamma_a) = \mathbb{E}[G'_n(\Gamma_a)] + \sum_{i=1}^{k_a} I_i(h_{a,n,i}),
\end{equation}
where $I_i$ indicates a multiple Wiener-It\^o integral of order $i$ with respect to the centered Poisson measure $\hat{\eta}_n$, $\E[G'_n(\Gamma_a)] = \int_{(\R^m)^{k_a}} h_{\Gamma_a,t_n}d\mu_n^{k_a}$, and, for $i=1,...,k_a$,  
\begin{equation}\label{defhani}
h_{a,n,i}(x_1,...,x_i) = \binom{k_a}{i} \int_{(\R^m)^{k_a-i}} h_{\Gamma_a,t_n} (x_1,...x_i,y_1,...,y_{k_a-i}) \mu_n^{k_a-i}(dy_1,...,dy_{k_a-i}).
\end{equation}   
Note that $h_{a,n,k_a} = h_{\Gamma_a, t_n}$. According to Theorem \ref{t:mainboundintro}, our proof is concluded if we can show that the six quantities appearing in formulae (\ref{e:alpha1})--(\ref{e:gamma2}) all converge to zero, as $n\to\infty$, at a rate of the order of $O\left( n^{ - \frac{1}{2(k-1)}} + n^{-\frac{k-k_0}{4(k-1)}}\right)$, whenever one selects the following arguments: (1) ${\bf F}_d = V'_n = (G'_n(\Gamma_1),...,G'_n(\Gamma_d))$, (2) ${\bf G}_m ={\bf G}_1 = \tilde{G}'_n(\Gamma_0)$, (3) $\lambda_i=\lambda_{n,i} = \E[G'_n(\Gamma_i)]$, $i=1,...,d$,  (4) $C=1$, and (5) $\mu = \mu_n$. We already know from \cite[Section 3]{LacPec_b} (see also Proposition \ref{p:g}-(a)) that the terms $\gamma_1(1, \tilde{G}'_n(\Gamma_0))$ and $\gamma_2(\tilde{G}'_n(\Gamma_0))$ both converge to zero at a rate $r_n$ such that
\[
r_n = O\left ( \frac{1}{\sqrt{n^{k_0}(t_n^m)^{k_0-1} }}\right).
\]
Since $ n^{k_0}(t_n^m)^{k_0-1} \sim n^{\frac{k-k_0}{k-1}}$, this implies that we only have to focus on the remaining four terms. We start by analysing the term $\alpha_1({\bf \lambda}_n, V'_n)$ and the first part of the term $\alpha_3(V'_n)$. 
\\~\\
Select $a,b = 1,...,d$. An application of the multiplication formula (\ref{e:product}), together with the definition of the derivative operator and the representation (\ref{e:proofc}), yields that 
\begin{eqnarray*}
&& \langle DG'_n(\Gamma_a), - DL^{-1}G'_n(\Gamma_b)\rangle_{L^2(\mu_n)} \\
&&= \sum_{i=1}^k \sum_{j=1}^k i \sum_{r=1}^{i\wedge j} (r-1)! \binom{i-1}{r-1}\binom{j-1}{r-1}\sum_{l=1}^r\binom{r-1}{l-1} I_{i+j-r-l}(h_{a,n,i}\star_r^l h_{b,n,j}) \\
&& = {\bf 1}_{\{a=b\}} \E[G'_n(\Gamma_a)] \\
&& \quad+ \sum_{i=1}^k \sum_{j=1}^k i \sum_{r=1}^{i\wedge j} (r-1)! \binom{i-1}{r-1}\binom{j-1}{r-1}\sum_{l=1}^r\binom{r-1}{l-1}{\bf 1}_{\{(i,j,r,l)\neq (k,k,k,k)\}} I_{i+j-r-l}(h_{a,n,i}\star_r^l h_{b,n,j}).
\end{eqnarray*}
Applying repeatedly the Cauchy-Schwarz inequality, one sees that, in order to prove that $\alpha_1({\bf \lambda}_n, V'_n)$ and $\alpha_3(V'_n)$ both converge to zero at the correct rate, it is sufficient to show that, for every $a,b=1,...,d$ and for every quadruple $(i,j,r,l)$ involved in the previous sum,
$$
\| h_{a,n,i}\star_r^l h_{b,n,j}\|_{L^2(\mu_n^{i+j-r-l})} = O\left(\sqrt{nt_n^m}\right) = O\left( n^{ - \frac{1}{2(k-1)}}\right)
$$ 
(the last equality is trivial). For any such $(i,j,r,l)$ we define the function $h_{a,b, t_n}^{(i,j,r,l)} : (\R^m)^{\alpha} \to \R$, where $\alpha = \alpha(i,j,r,l) = 4k-i-j-r+l$, as follows:
\begin{eqnarray}
h_{a,b, t_n}^{(i,j,r,l)}(x_1,...,x_\alpha) &=& h_{\Gamma_a,t_n}({\bf x}^{(1)}_{k-i}, {\bf x}^{(2)}_{i-r}, {\bf x}^{(3)}_{r-l}, {\bf x}^{(4)}_{l} )h_{\Gamma_b,t_n}({\bf x}^{(5)}_{k-j}, {\bf x}^{(6)}_{j-r}, {\bf x}^{(3)}_{r-l}, {\bf x}^{(4)}_{l} ) \\
&& \times h_{\Gamma_a,t_n}({\bf x}^{(7)}_{k-i}, {\bf x}^{(2)}_{i-r}, {\bf x}^{(3)}_{r-l}, {\bf x}^{(8)}_{l} )h_{\Gamma_b,t_n}({\bf x}^{(9)}_{k-j}, {\bf x}^{(6)}_{j-r}, {\bf x}^{(3)}_{r-l}, {\bf x}^{(8)}_{l} ),
\end{eqnarray}
where the bold letters represent multidimensional variables providing a lexicographic decomposition of $(x_1,...,x_\alpha)$. For instance, one has that ${\bf x}^{(1)}_{k-i} = (x_1,...,x_{k-i})$, ${\bf x}^{(2)}_{i-r} = (x_{k-i+1},...,x_{k-r})$, and so on, in such a way that $({\bf x}^{(1)}_{k-i}, {\bf x}^{(2)}_{i-r}, {\bf x}^{(3)}_{r-l}, {\bf x}^{(4)}_{l}, {\bf x}^{(5)}_{k-j}, {\bf x}^{(6)}_{j-r} ,{\bf x}^{(7)}_{k-i}, {\bf x}^{(8)}_{l} , {\bf x}^{(9)}_{k-j} )= (x_1,...,x_\alpha) $, and we set ${\bf x}^{(a)}_{p}$ equal to the empty set whenever $p=0$. Observe that each function $h_{a,b, t_n}^{(i,j,r,l)}$ is bounded by $1/k!^4$, and that the connectedness of the graphs $\Gamma_a, \Gamma_b$ yields that the mapping $(x_2,...,x_\alpha) \mapsto h_{a,b, 1}^{(i,j,r,l)}(0,x_2,...,x_\alpha)$, where $0$ stands for the origin, has compact support. Writing explicitly the squared contractions inside the integral, one sees that $\| h_{a,n,i}\star_r^l h_{b,n,j}\|^2_{L^2(\mu_n^{i+j-r-l})}$ is a multiple (with coefficient independent of $n$) of
\[
n^{\alpha}\int_{(\R^m)^{\alpha}} h_{a,b, t_n}^{(i,j,r,l)}(x_1,...,x_\alpha) f(x_1)\cdots f(x_\alpha)dx_1\cdots dx_\alpha.
\] 
Applying the change of variables $x_1 = x$ and $x_i = t_ny_i +x$, for $i=2,...,\alpha$, the above expression becomes
\begin{eqnarray*}
&& n^{\alpha}(t_n^m)^{\alpha-1} \int_{\R^m}f(x)  \int_{(\R^m)^{\alpha -1}} h_{a,b, 1}^{(i,j,r,l)}(0,y_2,...,y_\alpha) f(x+t_n y_2)\cdots f(x+t_ny_\alpha)dxdy_2\cdots dy_\alpha.
\end{eqnarray*}
Since, by dominated convergence, the integral on the RHS in the previous equation converges to the constant
\[
\int_{\R^m}f^\alpha (x)dx  \int_{(\R^m)^{\alpha -1}} h_{a,b, 1}^{(i,j,r,l)}(0,y_2,...,y_\alpha)dy_2\cdots dy_\alpha,
\]
we deduce that $\| h_{a,n,i}\star_r^l h_{b,n,j}\|^2_{L^2(\mu_n^{i+j-r-l})} = O\left(  n^{\alpha}(t_n^m)^{\alpha-1}\right) $. Since $$n^{\alpha}(t_n^m)^{\alpha-1} = n^k(t_n^m)^{k-1}(nt_n^m)^{\alpha-k}$$ and $\alpha-k\geq 1$ for every possible choice of $i,j,r,l$, we immediately deduce the desired conclusion for $\alpha_1({\bf \lambda}_n, V'_n)$ and the first part of $\alpha_3(V'_n)$.
\\~\\
To deal with $\alpha_2(V'_n)$ and the second part of $\alpha_3(V'_n)$, we apply the Cauchy-Schwarz inequality to write, for every $a,b=1,...,d$,
\[
\mathbb{E} \int_{\R^m}\left\lvert D_{z}G'_n(\Gamma_a)\left(D_{z}G'_n(\Gamma_a) - 1 \right)D_{z}L^{-1}G'_n(\Gamma_b) \right\rvert \mu_n(dz) \leq \sqrt{A(a,n) \times B(b,n)},
\]
where $A(a,n) = \mathbb{E} \int_{\R^m}  D_{z}G'_n(\Gamma_a)^2\left(D_{z}G'_n(\Gamma_a) - 1 \right)^2  \mu_n(dz)$ and $$B(b,n) =\mathbb{E} \int_{\R^m}[ D_{z}L^{-1}G'_n(\Gamma_b) ]^2 \mu_n(dz).$$ One can easily verify that the sequence $n\mapsto B(b,n)$ is bounded (this is a consequence of the fact that $G'_n(\Gamma_b)$ lives in a finite sum of Wiener chaoses). It follows that, in order to obtain the desired rate of convergence for this part, we just have to prove that, as $n\to\infty$, $A(a,n) = O\left( nt_n^m\right) $. To do this, one applies again the multiplication formula (\ref{e:product}) (for every fixed $z\in \R^m$) to deduce that, by virtue of the fact that $h_{\Gamma_a,t_n}$ has the special form of an indicator multiplied by the factor $k!^{-1}$, 
\begin{eqnarray}\label{e:boo}
&& D_z G'_n(\Gamma_a)(D_z G'_n(\Gamma_a)-1)   \\
&& = \sum_{i=1}^k \sum_{j=1}^k ij\!\! \sum_{r=0}^{i\wedge j -1} \!\!r! \binom{i-1}{r}\binom{j-1}{r}\sum_{l=0}^r\binom{r}{l}{\bf 1}_{\{(i,j,r,l)\neq (k,k,k-1,0)\}} \times \\
&& \quad\quad\quad \times I_{i+j-2-r-l}(h_{a,n,i}(z,\cdot) \star_r^l h_{a,n,j}(z,\cdot))\notag \\ 
&&\quad\quad\quad\quad\quad\quad\quad\quad\quad\quad\quad\quad - \sum_{t=1}^{k-1} tI_{t-1} (h_{a,n,t}(z,\cdot)) := \sum_{\gamma\in U} \xi_\gamma(z). \notag
\end{eqnarray}
In the last equality, the set $U$ represents the class of all indices $(i,j,r,l)$ and $t$ involved in the representation of $D_z G'_n(\Gamma_a)(D_z G'_n(\Gamma_a)-1)$, whereas $\xi_\gamma$ is the corresponding multiple integral process multiplied by the appropriate coefficient. To conclude, we apply the triangle inequality to deduce that 
\begin{equation}
\label{representationAan}
\sqrt{A(a,n)}\leq\sum_{\gamma\in U} \sqrt{\E\left[ \int_{\R^m} \xi_\gamma^2 (z) \m_n(dz)\right]}.
\end{equation}
We will show how to deal with the quadruple $(i,j,r,l) = (k,k,k-1,k-1)$, which requires additional arguments than the others (which can be addressed in a straightforward way). In the particular case where $(i,j,r,l) = (k,k,k-1,k-1)$, we are looking at the term $$\xi_{k,k,k-1,k-1} (z) = k^2 (k-1)! h_{a,n,k}(z,\cdot)\star_{k-1}^{k-1} h_{a,n,k}(z,\cdot).$$ Thus, we have 
\begin{eqnarray*}
\E\left[ \int_{\R^m} \xi_{k,k,k-1,k-1}^{2} (z) \m_n(dz)\right] &=& kk!\int_{\R^m}\left[h_{a,n,k}(z,\cdot)\star_{k-1}^{k-1} h_{a,n,k}(z,\cdot) \right]^{2} \mu_n(dz) \\
&=& kk!\int_{\R^m}\left(\int_{\left( \R^m\right)^{k-1} }h_{a,n,k}^{2}(z,y_1,...,y_{k-1})\mu_n^{k-1}(dy_1,...,dy_{k-1}) \right)^{2}  \mu_n(dz).
\end{eqnarray*}
Using the fact that $h_{a,n,k} = h_{\Gamma_a, t_n}$ along with the fact that $h_{\Gamma_a, t_n}$ has the form of an indicator function and finally recalling the definition of $ h_{a,n,1}$ given by (\ref{defhani}), we can write
\begin{equation*}
\E\left[ \int_{\R^m} \xi_{k,k,k-1,k-1}^2 (z) \m_n(dz)\right] = (k-1)!\int_{\R^m}h_{a,n,1}^{2}(z)\mu_n(dz) = (k-1)!\Vert h_{a,n,1} \Vert_{L^{2}(\mu_n)}^{2}.
\end{equation*}
The analysis of the contraction carried out in the previous steps of the proof allow us to conclude that this quantity goes to zero at the correct rate when $n$ goes to infinity (it corresponds to the $(1,1,1,1)$--contraction). Representing each remaining expectations in (\ref{representationAan}) as a contraction, and applying a change of variables analogous to the one described above gives the global and desired rate of convergence for $\alpha_2(V'_n)$ as well as for the second part of $\alpha_3(V'_n)$.
\\~\\
We now deal with the third and last part of $\alpha_3(V'_n)$. Applying the same strategy, we can write, for every $a,b,c=1,...,d$ with $a \neq b$,
\[
\mathbb{E} \int_{\R^m}\left\lvert D_{z}G'_n(\Gamma_a)D_{z}G'_n(\Gamma_b)D_{z}L^{-1}G'_n(\Gamma_c) \right\rvert \mu_n(dz) \leq \sqrt{C(a,b,n) \times D(c,n)},
\]
where $C(a,b,n) = \mathbb{E} \int_{\R^m}  D_{z}G'_n(\Gamma_a)^2 D_{z}G'_n(\Gamma_b)^2  \mu_n(dz)$ and $$D(c,n) =\mathbb{E} \int_{\R^m}[ D_{z}L^{-1}G'_n(\Gamma_c) ]^2 \mu_n(dz).$$ Again, the sequence $n\mapsto D(c,n)$ is bounded and we can write, for $a \neq b$,
\begin{eqnarray}
\label{e:boo2}
&& D_z G'_n(\Gamma_a)D_z G'_n(\Gamma_b)  \\
&& = \sum_{i=1}^k \sum_{j=1}^k ij\!\! \sum_{r=1}^{i\wedge j} (r-1)! \binom{i-1}{r-1}\binom{j-1}{r-1}\sum_{l=1}^r\binom{r-1}{l-1}I_{i+j-r-l}(h_{a,n,i}(z,\cdot) \star_{r-1}^{l-1} h_{b,n,j}(z,\cdot)) \\
&& := \sum_{\gamma\in I} \zeta_\gamma(z). \notag
\end{eqnarray}
In the last equality, the set $I$ represents the class of all indices $(i,j,r,l)$ involved in the representation of $D_z G'_n(\Gamma_a)D_z G'_n(\Gamma_b)$, whereas $\zeta_\gamma$ is the corresponding multiple integral process multiplied by the appropriate coefficient. This case is very similar to the previous one and the techniques used to prove that each of these expectation converge to zero as the correct rate are the same. However, there is one additional term that was not present in the case of $\alpha_2(V'_n)$. This is the term corresponding to the quadruple $(i,j,r,l) = (k,k,k-1,0)$. We will detail this particular case. We have
\begin{eqnarray*}
&& \E\left[ \int_{\R^m} \zeta_{k,k,k-1,0}^{2} (z) \m_n(dz)\right] = \E\int_{\R^m}(kk!)^2 I_{k-1}^{2}\left(h_{a,n,k}(z,\cdot) \star_{k-1}^{0} h_{b,n,k}(z,\cdot) \right)  \mu_n(dz) \\
&& = kk!^3\int_{\R^m}\int_{(\R^m)^{k-1}}h_{\Gamma_a, t_n}^{2}(z,y_1,...,y_{k-1})h_{\Gamma_b, t_n}^{2}(z,y_1,...,y_{k-1})  \mu_n^{k-1}(dy_1,...,dy_{k-1})\mu_n(dz).
\end{eqnarray*}
Using the fact that $h_{\Gamma_a, t_n}$ and $h_{\Gamma_b, t_n}$ have the form of indicator functions, we finally get
\begin{equation*}
\E\left[ \int_{\R^m} \zeta_{k,k,k-1,0}^{2} (z) \m_n(dz)\right] = kk!^{3} \left\langle h_{\Gamma_a, t_n},h_{\Gamma_b, t_n} \right\rangle _{L^{2}(\mu_n^{k})},
\end{equation*}
which is zero because $\Gamma_a$ and $\Gamma_b$ are not isomorphic ($h_{\Gamma_a, t_n}$ and $h_{\Gamma_b, t_n}$ cannot be non--zero at the same time or $\Gamma_a$ and $\Gamma_b$ would both be isomorphic to the same graph rendering them isomorphic to one another). This concludes the analysis of the term $\alpha_3(V'_n)$.
\\~\\
It remains to deal with $\beta(V'_n, \tilde{G}'_n(\Gamma_0))$. Using relation (\ref{e:uff}) with $\epsilon = 3$, we can write
$$\E \left\langle \vert DG'_n(\Gamma_a) \vert, \vert DL^{-1}\tilde{G}'_n(\Gamma_0) \vert \right\rangle _{L^{2}(\mu_{n})} \leq \E \left[\left[ DG'_n(\Gamma_a)\right] ^{2} \right]^{\frac{3}{4}} \times \E \left[\left[ DL^{-1}\tilde{G}'_n(\Gamma_0)\right] ^{4} \right]^{\frac{1}{4}}.$$ The term $\E \left[\left[ DG'_n(\Gamma_a)\right] ^{2} \right]^{\frac{3}{4}}$ is bounded and it remains to show that the term $\E \left[\left[ DL^{-1}\tilde{G}'_n(\Gamma_0)\right] ^{4} \right]^{\frac{1}{4}}$ goes to zero as $n$ goes to infinity. For this, we will refer to \cite[Section 3]{LacPec_b} where the rate of convergence of the term $\gamma_2(\tilde{G}'_n(\Gamma_0))$ is obtained by bounding it by a constant multiplied by $\sqrt{\E\int_{\R^m}\left[ DL^{-1}\tilde{G}'_n(\Gamma_0)\right] ^{4}\mu_{n}(dz)}$. It is then showed that this last term goes to zero at a rate of $O\left(n^{-\frac{k-k_0}{2(k-1)}} \right)$. The difference here lies in the fact that the square root is replaced by a power $\frac{1}{4}$, yielding a rate of convergence of $O\left(n^{-\frac{k-k_0}{4(k-1)}} \right)$.
\\~\\
When putting together all the rates of convergence for the different terms in the general bound, one sees that 
\begin{equation*}
d_{\star}(V_n, H_n) \leq A \sqrt{n t_n^m} + B n^{-\frac{k-k_0}{4(k-1)}} = O\left( n^{ - \frac{1}{2(k-1)}} + n^{-\frac{k-k_0}{4(k-1)}}\right) ,
\end{equation*}
where $A$ and $B$ are positive constants that do not depend on $n$. This concludes the proof.\qed

\medskip

\noindent{\bf Acknowledgments}. We thank Christoph Th\" ale for useful discussions.

\begin{alphasection}

\section{Appendix}

Throughout the Appendix,  $(Z,\mathscr{Z})$ denotes a Borel space endowed with a non-atomic $\sigma$-finite Borel measure $\mu$. We write $\eta$ to indicate a Poisson measure on $Z$ with control $\mu$. As in the main text, $\eta$ is assumed to be defined on some probability space $(\Omega, \mathscr{F}, \mathbb{P})$ such that $\mathscr{F}$ is the $\mathbb{P}$-completion of $\sigma(\eta)$. We also write $L^2(\mathbb{P}) = L^2(\Omega, \mathscr{F},\mathbb{P})$.   

\subsection{Malliavin operators}\label{app:mall}
We now define some Malliavin-type operators associated with the Poisson measure $\eta$. We follow the work by Nualart and Vives \cite{nuaviv}.

\medskip

\noindent \underline{\bf The derivative operator $D$}. 

\smallskip

\noindent For every $F\in L^2(\mathbb{P})$, the derivative of $F$, $DF$ is defined as an element of $L^2(\mathbb{P};L^2(\mu))$, that is, of the space of the jointly measurable random functions $u:\Omega \times Z \to \mathbb{R}$ such that $E \left[\int_Z u_z^2 \mu(dz) \right] <\infty$.

\begin{definition}{\rm 
\begin{enumerate}
 \item The domain of the derivative operator $D$, written  ${\rm dom} D$, is the set of all random variables $F\in L^2(P)$ admitting a chaotic decomposition \eqref{e:chaos} such that
$$ \sum_{k\geq 1} k k!\|f_k \|^2_{L^2(\mu^k)} < \infty ,$$
 \item For any $F\in {\rm dom}D$, the random function $z \mapsto D_z F$ is defined by
$$ D_z F= \sum_{k \geq 1}^{\infty} k I_{k-1}(f_k(z,\cdot)) .$$
 \end{enumerate}}
\end{definition}

\medskip

\noindent \underline{\bf The divergence operator $\delta$}.  

\smallskip

\noindent Thanks to the chaotic representation property of $\eta$, every random function
$u \in L^2(\mathbb{P},L^2(\mu))$ admits a unique representation of the type
\begin{equation} \label{skor}
 u_z = \sum_{k \geq 0}^{\infty}  I_{k}(f_k(z,\cdot)) ,\,\, z\in Z,
\end{equation}
where the kernel $f_k$ is a function of $k+1$ variables, and $f_k(z,\cdot)$ is an element of $L^2_s(\mu^k)$. The {\sl divergence operator} $\delta(u)$ maps a random function $u$ in its domain to an element of $L^2(P)$.\\

\begin{definition}{\rm
\begin{enumerate}
  \item  The domain of the divergence operator, denoted by  ${\rm dom} \delta$, is the collection of all $u\in L^2(P,L^2(\mu))$ having the above chaotic expansion (\ref{skor}) satisfied the condition:
$$ \sum_{k\geq 0}  (k+1)! \|f_k \|^2_{L^2(\mu^(k+1))} < \infty. $$
  \item For $u\in {\rm dom}\delta$, the random variable $\delta(u)$ is given by
      $$ \delta (u) = \sum_{k\geq 0} I_{k+1}(\tilde{f}_k), $$
      where $\tilde{f}_k$ is the canonical symmetrization of the $k+1$ variables function $f_k$.
\end{enumerate}}
\end{definition}
As made clear in the following statement, the operator $\delta$ is indeed the adjoint operator of $D$.
\begin{lemma}[Integration by parts]\label{L : IBP}
 For every $G\in {\rm dom} D$ and $u\in {\rm dom} \delta$, one has that
$$ \mathbb{E}[G \delta(u)] = \mathbb{E}[\langle D G, u \rangle_{L^2(\mu)}]. $$
\end{lemma}
The proof of Lemma \ref{L : IBP} is detailed e.g. in \cite{nuaviv}.\\

\medskip

\noindent \underline{\bf The Ornstein-Uhlenbeck generator $L$}.

\smallskip

\begin{definition}{\rm
\begin{enumerate}
  \item  The domain of the Ornstein-Uhlenbeck generator, denoted by  ${\rm dom} L$, is the collection of all $F \in L^2(\mathbb{P})$ whose chaotic representation \label{chaos} satisfies the condition:
$$ \sum_{k\geq 1}  k^2 k! \|f_k \|^2_{L^2(\mu^k)} < \infty $$
  \item The Ornstein-Uhlenbeck generator $L$ acts on random variable $F\in {\rm dom}L$ as follows:
      $$ LF = - \sum_{k\geq 1} k I_{k}(f_k) .$$
\end{enumerate}}
\end{definition}

\medskip

\noindent \underline{\bf The pseudo-inverse of $L$}.

\smallskip

\begin{definition}{\rm
\begin{enumerate}
  \item  The domain of the pseudo-inverse of the Ornstein-Uhlenbeck generator, denoted by $L^{-1}$, is the space $L^2_0(\mathbb{P})$ of \it{centered} random variables in $L^2(\mathbb{P})$.
  \item For $F = \sum\limits_{k\geq 1} I_k (f_k) \in L^2_0(\mathbb{P})$ , we set
      $$ L^{-1}F = - \sum_{k\geq 1} \cfrac{1}{k} I_{k}(f_k). $$
\end{enumerate}}
\end{definition}

\subsection{Contractions}\label{app:cont}

{ Contraction operators} play a crucial role in multiplication formulae and in the computation of expectations involving powers of functionals of the Poisson measure $\eta$. In what follows, we shall define these operators and discuss some of their basic properties. The reader is referred e.g. to \cite[Sections 6.2 and 6.3]{PeTa} for further details.

\smallskip

\noindent The kernel $f \star_r^l g$ on $Z^{p+q-r-l}$, associated with functions $f\in L^2_s(\mu^p) $ and $g \in L^2_s(\mu^q) $, where $p,q \geq 1$, $r=1,\ldots, p\wedge q$ and $l=1,\ldots,r $, is defined as follows:
\begin{eqnarray}
& & f \star_r^l
g(\gamma_1,\ldots,\gamma_{r-l},t_1,\ldots,t_{p-r},s_1,\ldots,s_{q-r}) \label{contraction} \\
&=& \int_{Z^l} \mu^l(dz_1,\ldots,dz_l)
f(z_1,\ldots,z_l,\gamma_1,\ldots,\gamma_{r-l},t_1,\ldots,t_{p-r}) \nonumber \\
& & \quad\quad\quad\quad\quad\quad\quad\quad\quad\quad\quad\quad \times g(z_1,\ldots,z_l,\gamma_1,\ldots,\gamma_{r-l},s_1,\ldots,s_{q-r}). \nonumber
\end{eqnarray}
Roughly speaking, the star operator `$\,\star_r^l\,$' reduces the number of variables in the tensor product of $f$ and $g$ from $p+q$ to $p+q-r-l$: this operation is realized by first identifying $r$ variables in $f$ and $g$, and then by integrating out $l$ among them. To deal with the case $l=0$ for $r=0,\ldots, p\wedge q$, we set
\begin{eqnarray*}
& &f \star_r^0
g(\gamma_1,\ldots,\gamma_{r},t_1,\ldots,t_{p-r},s_1,\ldots,s_{q-r}) \\
&=& f(\gamma_1,\ldots,\gamma_{r},t_1,\ldots,t_{p-r})
g(\gamma_1,\ldots,\gamma_{r},s_1,\ldots,s_{q-r}),
\end{eqnarray*}
 and
$$ f \star_0^0 g (t_1,\ldots,t_{p},s_1,\ldots,s_{q}) =f\otimes g (t_1,\ldots,t_{p},s_1,\ldots,s_{q})= f(t_1,\ldots,t_{p})
g(s_1,\ldots,s_{q}). $$

\noindent The kernel $ f \star_r^l g$ is called the {\it contraction of index $(r,l)$} between $f$ and $g$. The above introduced `star notation' is standard,
and has been first used by Kabanov in \cite{Kab} (see also
Surgailis \cite{Surg1984}). Plainly, for some choice of $f,g,r,l$
the contraction $f\star^l_r g$ may not be well-defined. The contractions of the following three types
are well-defined (although possibly infinite) for every $1\leq p\leq q$
and every pair of kernels $g\in L_s^2(\mu^p), \, f\in L_s^2(\mu^q)$:
\begin{itemize}
\item[(a)] $f\star_r^0 g(z_1,\ldots,z_{p+q-r})$, where $r=0,\ldots,p$;
\item[(b)] $f\star_q^l f(z_1,\ldots,z_{q-l}) = \int_{Z^l} f^2(z_1,\ldots,z_{q-l},\cdot)d\mu^l $, for every
$l=1,\ldots,q$;
\item[(c)] $f\star_r^r g$, for $r=0,\ldots,p$.
\end{itemize}

\noindent In particular, a contraction of the type $f\star_q^l f$,
where $l=1,\ldots,q-1$ may equal $+\infty$ at some point
$(z_1,\ldots,z_{q-l})$. The following (elementary) statement ensures
that any kernel of the type $f\star_r^r g$ is square-integrable.

\begin{lemma}\label{L : star-integrability}
Let $p,q\geq 1$, and let $f\in L_s^2(\mu^q)$ and $g\in
L_s^2(\mu^p)$. Fix $r=0,\ldots,q\wedge p$. Then, $f\star_r^r g \in
L^2(\mu^{p+q-2r})$.
\end{lemma}

\end{alphasection}

\end{document}